\documentclass{amsart}
\usepackage{xspace, graphicx, epsfig, amssymb, amsmath, amsthm, stmaryrd}

\newcommand {\Z} {\mathbb{Z}}

\newcommand {\C} {\mathbb{C}}

\newcommand {\Proj} {\mathbb{P}}

\newcommand {\Aut} {\mathrm{Aut}}
\newcommand {\End} {\mathrm{End}\:}
\newcommand {\Hom} {\mathrm{Hom}}

\newcommand {\cod} {\mathrm{codim\ }}

\newcommand {\tr} {\mathrm{Tr}}

\newcommand{\sgn}{\mathrm{sgn}}

\newcommand {\GL} {\mathrm{GL}}

\newcommand {\SO} {\mathrm{SO}}
\newcommand {\SP} {\mathrm{SP}}

\def\a{\mathfrak {a}}
\def\b{\mathfrak {b}}
\def\gl{\mathfrak {gl}}
\def\h{\mathfrak {h}}
\def\g{\mathfrak {g}}
\def\l{\mathfrak {l}}
\def\m{\mathfrak {m}}

\def\sl{\mathfrak{sl}}
\def\so{\mathfrak{so}}
\def\sp{\mathfrak{sp}}
\def\t{\mathfrak{t}}
\def\z{\mathfrak {z}}
\def\k{\mathfrak{k}}
\def\s{\mathfrak{s}}

\newcommand {\ad}{\mathrm{ad \:}}
\newcommand {\Span}{\mathrm{Span}}
\newcommand {\supp}{\mathrm{supp \:}}
\newcommand{\Sym}{\mathrm{Sym}}
\newcommand {\rad}{\mathrm{rad \:}}
\newcommand {\rank}{\mathrm{rank \:}}

\newtheorem{thm}{Theorem}[section]
\newtheorem{lemma}[thm]{Lemma}
\newtheorem{prop}[thm]{Proposition}

\newtheorem{defn}[thm]{Definition}
\newtheorem{cor}[thm]{Corollary}

\begin{document}

\title{Cartan subalgebras of root-reductive Lie algebras}  

\author{Elizabeth Dan-Cohen}
\address{E. D.-C. Department of Mathematics \\ University of California at Berkeley \\ Berkeley, California 94720, USA}
\email{edc@math.berkeley.edu}

\author{Ivan Penkov} 
\address{I.P. International University Bremen\\ P.O. Box 750 561\\ D-28725 \\ Bremen, Germany}
\email{i.penkov@iu-bremen.de}

\author{Noah Snyder}
\address{N.S. Department of Mathematics \\ University of California at Berkeley \\ Berkeley, California 94720, USA}
\email{nsnyder@math.berkeley.edu}

\begin{abstract}
Root-reductive Lie algebras are direct limits of finite-dimensional reductive Lie algebras under injections which preserve the root spaces.  It is known that a root-reductive Lie algebra is a split extension of an abelian Lie algebra by a direct sum of copies of finite-dimensional simple Lie algebras as well as copies of the three simple infinite-dimensional root-reductive Lie algebras $\sl_\infty$, $\so_\infty$, and $\sp_\infty$.  As part of a structure theory program for root-reductive Lie algebras, Cartan subalgebras of the Lie algebra $\gl_\infty$ were introduced and studied in \cite{N-P}.

In the present paper we refine and extend the results of \cite{N-P} to the case of a general root-reductive Lie algebra $\g$.  We prove that the Cartan subalgebras of $\g$ are the centralizers of maximal toral subalgebras and that they are nilpotent and self-normalizing.  We also give an explicit description of all Cartan subalgebras of the simple Lie algebras $\sl_\infty$, $\so_\infty$, and $\sp_\infty$.

We conclude the paper with a characterization of the set of conjugacy classes of Cartan subalgebras of the Lie algebras $\gl_\infty$, $\sl_\infty$, $\so_\infty$, and $\sp_\infty$ with respect to the group of automorphisms of the natural representation which preserve the Lie algebra.
\end{abstract}

\maketitle

\section{Introduction} 
In this paper we continue the study of Cartan subalgebras of infinite-dimensional root-reductive Lie algebras initiated in \cite{N-P}.  We refine and extend the description of Cartan subalgebras of $\gl_\infty$ given in \cite{N-P} to the case of a general root-reductive Lie algebra.  We also solve the problem of describing the set of conjugacy classes of Cartan subalgebras of the three simple root-reductive Lie algebras $\sl_\infty$, $\so_\infty$, $\sp_\infty$, as well as of $\gl_\infty$, under the group of automorphisms of the natural representation which preserve the Lie algebra.  Root-reductive Lie algebras are a specific class of locally finite Lie algebras, and as background literature on locally finite Lie algebras we recommend \cite{BahturinBenkart}, \cite{BahturinStrade}, \cite{B1}, \cite{B2}, \cite{BaranovZhilinski}, and \cite{D-P1}.

Here is a description of the contents of the paper.  Its first main part, Sections \ref{sectiona} and \ref{sectionb}, is concerned with the definition and description of Cartan subalgebras.  According to \cite{N-P}, a Cartan subalgebra of $\gl_\infty$ can be defined as a locally nilpotent subalgebra which is the centralizer of the set of all semisimple parts of its own elements.  We accept this definition for any locally reductive Lie algebra $\g$ and prove that similarly to the case of $\gl_\infty$, Cartan subalgebras are nothing but centralizers of arbitrary maximal toral subalgebras.  In particular, a Cartan subalgebra is always self-normalizing.  We then prove that a Cartan subalgebra can be characterized equivalently as a subalgebra $\h \subseteq \g$ which coincides with the set of $\h_{fin}$-locally nilpotent vectors in $\g$, where $\h_{fin}$ runs over all finite-dimensional subalgebras of $\h$.  In addition, we show that if $\g$ is a root-reductive Lie algebra, any Cartan subalgebra $\h$ is nilpotent and coincides with the set of all $\h$-locally nilpotent vectors in $\g$.  The main new phenomenon in the case of a general root-reductive Lie algebra versus the case of $\gl_\infty$ is that Cartan subalgebras are no longer necessarily commutative and that the adjoint action of a Cartan subalgebra $\h$ on itself no longer has to be locally finite. 

We treat in detail the simple Lie algebras $\sl_\infty$, $\so_\infty$, and $\sp_\infty$.  The case of $\sl_\infty$ is very similar to the case of $\gl_\infty$ considered in \cite{N-P}.  The main new result in the other two cases is that maximal toral subalgebras are in one-to-one correspondence with maximal orthogonal or symplectic self-dual systems of $1$-dimensional subspaces in the natural representation.  Similarly to the case of $\gl_\infty$, such maximal self-dual systems do not have to span the natural representation, a phenomenon that gives rise to nonsplitting maximal toral subalgebras.  In the case of $\so_\infty$ the centralizer of a maximal toral subalgebra is in general nilpotent of depth no greater than $2$ but not necessarily commutative.  Remarkably enough, the ideal of nilpotent elements of a Cartan subalgebra of $\so_\infty$ is itself a nilpotent orthogonal Lie algebra of a vector space with a degenerate symmetric form, and all such nilpotent degenerate algebras occur inside Cartan subalgebras of $\so_\infty$.  In the case of $\sl_\infty$, the analogous nilpotent degenerate subalgebras do not occur inside Cartan subalgebras of $\sl_\infty$, and 
for $\gl_\infty$ and $\sp_\infty$, there are no analogous nilpotent nonabelian subalgebras, which is consistent with the fact that $\gl_\infty$, $\sl_\infty$, and $\sp_\infty$ admit only abelian Cartan subalgebras.

In Section \ref{sectionc} we address the conjugacy problem for Cartan subalgebras posed in \cite{N-P}.  More precisely, in the case of $\gl_\infty$ it was established that certain standard discrete invariants of Cartan subalgebras are not sufficient to characterize the conjugacy classes of Cartan subalgebras of $\gl_\infty$ with respect to the group of automorphisms of the natural representation which preserve $\gl_\infty$.  In this paper we describe the missing continuous invariants in terms of a linear algebraic structure which we call a complement datum.  Our main theorem in the second part, Theorem \ref{equivthm}, solves the conjugacy problem in terms of the combinatorics of complement data.  It turns out that for a generic set of fixed standard invariants, there are uncountably many conjugacy classes of Cartan subalgebras.  In all remaining cases, there are only finitely many conjugacy classes, and in fact at most $2$.  We give criteria to distinguish the uncountable and finite cases and provide representatives in each finite case.

\section{Preliminaries}

All  vector spaces and Lie algebras are defined over the field of complex numbers $\C$.  

A Lie algebra $\g$ is {\it locally finite} (respectively {\it locally nilpotent}) if every finite subset of 
$\g$ is contained in a finite-dimensional (resp. nilpotent) subalgebra.  If $\g$ is at most countable dimensional, then being locally finite is equivalent to being isomorphic to a union $\bigcup_{i \in \Z_{>0}} \g_i$ of nested finite-dimensional Lie subalgebras $\g_i \subseteq \g_{i+1}$.  A $\g$-module $M$ is {\it locally finite} if each element $m \in M$ is contained in a finite-dimensional submodule, and $M$ is {\it locally nilpotent} if for any $m \in M$ there exists an $i \in \Z_{>0}$ with $\g^i \cdot m = 0$. Furthermore, we say that a $\g$-module $M$ is a {\it generalized weight module} if $M = \bigoplus_{\lambda \in \g^*} M^\lambda(\g)$, where 
$$ M^\lambda (\g):= \{ m \in M \colon \exists i \in \Z_{>0} \: s.t. \; \forall x \in \g, (x - \lambda (x) 1)^i \cdot m = 0\}. $$
 We define a generalized weight module $M$ to be a {\it weight module} if in addition, 
 $$ M^\lambda (\g) =  \{ m \in M \colon \forall x \in \g, \ x \cdot m = \lambda(x) m\} $$ 
 for each $\lambda \in \g^*$. The {\it support} in $\g^*$ of a module $M$ is the set $$ \supp M := \{ \alpha \in \g^* \colon M^\alpha(\g) \neq 0\}. $$ 

For an arbitrary Lie subalgebra $\h \subseteq \g$, we set $\overline{\g^0(\h)} = \bigcap_i \g^0(\h_i)$ where $\bigcup_i \h_i = \h$ is any exhaustion of $\h$ by finite-dimensional subalgebras $\h_i$.  In other words, $\overline{\g^0(\h)}$ is the subalgebra of all $\h_{fin}$-locally nilpotent vectors in $\g$, where $\h_{fin}$ runs over all finite-dimensional subalgebras of $\g$.

We call a Lie algebra \emph{locally reductive} if it is the union $\bigcup_{i \in \Z_{>0}} \g_i$ of nested finite-dimensional reductive Lie algebras $\g_i$ such that $\g_i$ is reductive in $\g_{i+1}$.  Recall that for finite-dimensional Lie algebras we have a notion of a Jordan decomposition, and in particular of semisimple elements and nilpotent elements.  Since these notions are preserved under injections $\g_i \hookrightarrow \g_{i+1}$ where $\g_i$ is reductive in $\g_{i+1}$, for any locally reductive Lie algebra we can talk about a \emph{Jordan decomposition} of an element.  In addition, note that any nilpotent element $x$ of a locally reductive Lie algebra $\g$ lies in $[\g,\g]$.  Indeed, $x \in \g_i$ for some $i$, and thus $x \in [\g_i,\g_i] \subseteq [\g,\g]$.

We call a subalgebra $\k$ of a locally reductive Lie algebra $\g$ \emph{splittable} if for every $k \in \k$, both the semisimple and nilpotent Jordan components $k_{ss}$ and $k_{nil}$ of $k$ belong to $\k$.  A subalgebra $\t \subseteq \g$ is \emph{toral} if every element is semisimple.  Every toral subalgebra is abelian: the standard proof of this fact extends from the case of a reductive Lie algebra to a locally reductive Lie algebra, cf. \cite{N-P}.  By definition, a toral subalgebra $\t \subseteq \g$ is \emph{splitting} if $\g$ is a weight module with respect to the adjoint action of  $\t$ on $\g$.

For any subset $\a \subseteq \g$ and any subalgebra $\k \subseteq \g$, we define the centralizer of $\a$ in $\k$, denoted  $\z_\k(\a)$, to be the set of elements of $\k$ which commute in $\g$ with all elements of $\a$.  The center of $\g$ is denoted $\z(\g)$.

Let $V$ and $V_*$ be vector spaces of countable dimension with a nondegenerate pairing $\langle \cdot , \cdot \rangle : V \times V_* \to \C$.  We define the Lie algebra $\gl (V,V_*)$ to be the Lie algebra associated to the associative algebra $V \otimes V_*$ ($V \otimes V_*$ is an associative algebra whose multiplication satisfies $(v \otimes w) \cdot (v' \otimes w') = \langle v' , w \rangle v \otimes w'$ for $v,v' \in V$, $w,w' \in V_*$).  The derived subalgebra of $\gl(V,V_*)$ has codimension $1$ and is denoted $\sl(V,V_*)$.  If $V$ is endowed with a nondegenerate symmetric (resp. antisymmetric) bilinear form, then one can choose $V_* = V$ and denote by $\so(V)$ (resp. $\sp(V)$) the Lie subalgebra $\bigwedge^2 V$ (resp. $\Sym^2 (V)$) of the Lie algebra associated to the associative algebra $V \otimes V$.

A result of Mackey \cite[p. 171]{Mackey} implies that all associative algebras $V \otimes V_*$ are isomorphic, as long as $\dim V = \dim V_* = \aleph_0$ and the pairing $V \times V_* \rightarrow \C$ is nondegenerate.  Hence, up to isomorphism, the Lie algebras $\gl(V,V_*)$, $\sl(V,V_*)$, $\sl(V)$, and $\sp(V)$ do not depend on $V$, and the usual representatives of these isomorphism classes are denoted $\gl_\infty$, $\sl_\infty$, $\so_\infty$, and $\sp_\infty$, respectively.  Clearly the Lie algebras $\gl_\infty$, $\sl_\infty$, $\so_\infty$, and $\sp_\infty$ are locally reductive, as one can also choose exhaustions $V = \bigcup_{i \in Z_{>0}} V_i$ and $V_* = \bigcup_{i \in Z_{>0}} W_i$ (where $V_i = W_i$ for $V = V_*$) such that the pairings $V_i \times W_i \rightarrow \C$ are nondegenerate.

\section{Cartan subalgebras of locally reductive Lie algebras} \label{sectiona}

Let $\g$ be a locally reductive Lie algebra.  If $\h$ is a subalgebra of $\g$, let $\h_{ss}$ denote the set of semisimple Jordan components of the elements of $\h$.  Following \cite{N-P}, we say that a subalgebra $\h$ of $\g$ is a \emph{Cartan subalgebra} if $\h$ is locally nilpotent and $\h = \z_\g(\h_{ss})$.  

We start with the following proposition which generalizes a result of \cite{N-P}. 

\begin{prop} \label{centralizer}
Let $\h$ be a locally nilpotent subalgebra of a locally reductive Lie
algebra $\g$.  Then the following assertions hold:
\begin{enumerate}
\item $\h \subseteq \z_\g(\h_{ss})$;
\item \label{toral} $\h_{ss}$ is a toral subalgebra of $\g$;
\item  \label{self-normalizing} $\z_\g (\h_{ss})$ is a self-normalizing subalgebra of $\g$.
\end{enumerate}
\end{prop}

\begin{proof}
Let $h, h' \in \h$.  The local nilpotence of $\h$ implies that $(\ad
h)^n(h') = 0$ for some $n$.  Since $\ad h_{ss}$ is a polynomial with no constant term in $\ad h$, it follows that $(\ad h_{ss})(\ad h)^{n-1}(h') = 0$.  Because an element commutes with its semisimple part, $(\ad h)^{n-1}(\ad h_{ss})(h') = 0$, and it follows by induction that $(\ad h_{ss})^n(h') = 0$.  Hence $(\ad h_{ss})(h') = 0$.  Thus $\h \subseteq \z_\g(\h_{ss})$.

Furthermore, by the same argument, $(\ad h')(h_{ss}) = 0$ implies that $(\ad h'_{ss})(h_{ss}) = 0$.  Therefore any two elements of $\h_{ss}$ commute.  Since the sum of any two commuting semisimple elements is semisimple, $\h_{ss}$ is a subalgebra.

Finally, suppose $x$ is in the normalizer of $ \z_\g(\h_{ss})$.  For any $y \in \h_{ss}$, we have that $[x,y] \in \z_\g(\h_{ss})$.  Thus $[[x,y],y]=0$, and as $y$ is semisimple it follows that $[x,y]=0$. Hence $x \in \z_\g(\h_{ss})$, i.e. $\z_\g(\h_{ss})$ is self-normalizing.
\end{proof}

The following theorem is our main general result characterizing Cartan subalgebras of locally reductive Lie algebras.  It generalizes a well-known result for finite-dimensional reductive algebras (cf. \cite[Ch. VII, \S 4]{Bourbaki}), as well as a result from \cite{N-P}.

\begin{thm} \label{main}
Let $\g$ be a locally reductive Lie algebra and $\h$ a subalgebra of $\g$.  The following conditions on $\h$ are equivalent:
\begin{enumerate}
\item \label{C} $\h$ is a Cartan subalgebra;
\item \label{D} $\h = \z_\g(\h_{ss})$ and $\h_{ss}$ is a subalgebra;
\item  \label{E} $\h = \z_\g(\t)$ for some maximal toral subalgebra $\t \subseteq \g$;
\item \label{F} $\h = \overline{\g^0(\h)}$.
\end{enumerate}
In addition, any Cartan subalgebra $\h$ is splittable and self-normalizing.
\end{thm}

We first prove a short lemma.

\begin{lemma}\label{lemma1}
If $\h$ is locally nilpotent and splittable, then $\overline{\g^0(\h)} = \z_\g(\h_{ss})$.
\end{lemma}
\begin{proof}
By \cite[Ch. VII, \S 5, Prop. 5]{Bourbaki} $\h = \h_{ss} \oplus \h_{nil}$, with $\h_{nil}$ being the subalgebra of all nilpotent elements in $\h$.  It follows that  $\overline{\g^0(\h)} = \overline{\g^0(\h_{ss})} \cap \overline{\g^0(\h_{nil})}$.  Since elements of $\h_{ss}$ are semisimple, $\overline{\g^0(\h_{ss})} = \z_\g(\h_{ss})$.  Clearly $\overline{\g^0(\h_{nil})} = \g$. Hence $\overline{\g^0(\h)} = \z_\g(\h_{ss})$.
\end{proof}

\begin{proof}[Proof of Theorem \ref{main}] Fix an exhaustion $\g = \bigcup_{i \in \Z_{>0}} \g_i$, where $\g_i$ is a finite-dimensional reductive Lie algebra with $\g_i$ reductive in $\g_{i+1}$. 

Clearly (\ref{C}) implies (\ref{D}), by Proposition \ref{centralizer} (\ref{toral}).  To show that (\ref{D}) implies (\ref{C}), we first prove that (\ref{D}) implies that $\h$ is splittable.  Suppose that $\h$ satisfies (\ref{D}).  For any $i \in \Z_{>0}$ note that $$\h \cap \g_i =  \z_\g(\h_{ss}) \cap \g_i =  \bigcap_{k \geq i}(\z_{\g_k}(\h_{ss} \cap \g_k) \cap \g_i).$$  Since $\dim \g_i < \infty$, we have
$\h \cap \g_i = \z_{\g_j}(\h_{ss} \cap \g_j) \cap \g_i$
for some sufficiently large $j \geq i$.
Since $\h_{ss} \cap \g_j$ is a subalgebra of $\g_j$, we know from \cite[Ch. VII, \S 5, Prop. 3 Cor. 1]{Bourbaki} that $\z_{\g_j}(\h_{ss} \cap \g_j)$ is a splittable subalgebra of $\g_j$.  The reductive Lie algebra $\g_i$ is also splittable, by \cite[Ch. VII, \S 5, Prop. 2]{Bourbaki}.  Then the intersection $\z_{\g_j}(\h_{ss} \cap \g_j) \cap \g_i$ is splittable, too.  Being a union of splittable algebras, $\h$ is splittable.  

To show the local nilpotence of $\h$, notice that the equality $\h = \z_\g(\h_{ss})$ implies that every element of $\h_{ss}$ is in the center of $\h$.  Now consider a general element $h = h_{ss} + h_{nil} \in \h$.  Choose $k$ such that $(\ad h_{nil})^k = 0$.  For any $x \in \h$, $$(\ad h)^k (x) = (\ad(h_{ss} + h_{nil}))^k(x) = (\ad h_{nil})^k(x) = 0.$$  Hence $\h$ is locally nilpotent, i.e. (\ref{D}) implies (\ref{C}).

Next, it is clear that (\ref{D}) implies (\ref{E}).  The equality $\h = \z_{\g}(\h_{ss})$ shows that any semisimple element of $\g$ which centralizes $\h_{ss}$ is already in $\h_{ss}$.  Thus $\h_{ss}$ is a maximal toral subalgebra and (\ref{E}) holds.

Let $\h$ satisfy (\ref{E}).  The same argument as above shows that $\h$ is splittable, whence $\h_{ss} \subseteq \h$.  Then clearly $\t \subseteq \h_{ss}$.  If $\t \neq \h_{ss}$, the existence of a semisimple element $h \in \h \setminus \t$ contradicts the maximality of $\t$.  Therefore $\t = \h_{ss}$, and (\ref{E}) implies (\ref{D}).

Note that (\ref{C}) implies (\ref{F}).  Indeed, let $\h$ be a Cartan subalgebra.  We know $\h$ is splittable and locally nilpotent, so by Lemma \ref{lemma1}, $\h = \z_\g(\h_{ss}) = \overline{\g^0(\h)}$.

To show that (\ref{F}) implies (\ref{C}), assume that $\h = \overline{\g^0(\h)}$.  Then clearly $\h$ is locally nilpotent.  An argument similar to that above shows that $\h$ is also splittable.  Indeed, for any $i \in \Z_{>0}$ we have 
$$\overline{\g^0(\h)} \cap \g_i = \bigcap_{k \geq i} \left(\g_k^0(\h \cap \g_k) \cap \g_i \right).$$  The finite dimensionality of $\g_i$ yields $\overline{\g^0(\h)} \cap \g_i = \g_{j}^0(\h \cap \g_{j}) \cap \g_i$ for some sufficiently large $j \geq i$.  It is well known that $\g_j^0(\h \cap \g_j)$ is a splittable subalgebra of $\g_j$, see \cite[Ch. VII, \S 1, Prop. 11]{Bourbaki}.  Since $\g_i$ is also splittable, the intersection $\g_{j}^0(\h \cap \g_{j}) \cap \g_i$ is splittable, too.   Hence $\overline{\g^0(\h)} \cap \g_i$ is splittable.  Being a union of splittable algebras, $\h$ is splittable.  Therefore Lemma \ref{lemma1} implies that $\h = \overline{\g^0(\h)} = \z_\g(\h_{ss})$.

In addition, by Proposition \ref{centralizer} (\ref{self-normalizing}), a subalgebra $\h$ satisfying (\ref{C}) is self-normalizing.  As we have already seen that a Cartan subalgebra is splittable, the proof of Theorem \ref{main} is complete.
\end{proof}

The most common characterization of a Cartan subalgebra of a finite-dimensional Lie algebra is as a nilpotent self-normalizing subalgebra.  According to Theorem \ref{main}, any Cartan subalgebra of a locally reductive Lie algebra is self-normalizing.  Any Cartan subalgebra of a root-reductive Lie algebra is in addition nilpotent.  In fact, we show in Theorem \ref{depth2} that it is nilpotent of depth at most $2$.  In contrast with the case of $\gl_\infty$ or $\sl_\infty$, where all Cartan subalgebras are abelian and hence act locally finitely on themselves for a trivial reason (cf. \cite{N-P}), the Cartan subalgebras of $\so_\infty$ do not necessarily act locally finitely on themselves, as discussed at the end of subsection \ref{sosubsection} below.  We do not know whether all nilpotent self-normalizing subalgebras of a root-reductive Lie algebra are Cartan subalgebras.  We also do not know whether a Cartan subalgebra of a general locally reductive Lie algebra is necessarily nilpotent.

In what follows we call a Cartan subalgebra $\h$ \emph{splitting} if $\h_{ss}$ is a splitting toral subalgebra of $\g$.

\section{Cartan subalgebras of root-reductive Lie algebras}\label{sectionb}

This section begins with some brief introductory material on root-reductive Lie algebras.  We then consider separately the cases of a root-reductive Lie algebra $\g$ containing a simple component of type $\sl_\infty$, $\so_\infty$, and $\sp_\infty$.  In particular, we describe all Cartan subalgebras of the simple Lie algebras $\sl_\infty$, $\so_\infty$, and $\sp_\infty$.  The section concludes with a general theorem about the Cartan subalgebras of an arbitrary root-reductive Lie algebra and with examples.

\subsection{Root-reductive Lie algebras}

An inclusion of finite-dimensional reductive Lie algebras $\l \subseteq \m$ is a \emph{root inclusion} if, for some Cartan subalgebra $\h_\m$ of $\m$, the subalgebra $\l \cap \h_\m$ is a Cartan subalgebra of $\l$ and each $\l \cap \h_\m$-root space $\l^\alpha$ is also a root space of $\m$.  Informally, root spaces of $\l$ are root spaces of $\m$.  A Lie algebra $\g$ is called \emph{root-reductive} if it is isomorphic to a union $\bigcup_{i \in \Z_{>0}} \g_i$ of nested reductive Lie algebras with respect to root inclusions for a fixed choice of nested Cartan subalgebras $\h_i \subseteq \g_i$ with $\h_{i-1} = \h_i \cap \g_{i-1}$.

Finite-dimensional reductive Lie algebras are root-reductive, as are $\gl_\infty$, $\sl_\infty$, $\so_\infty$, and $\sp_\infty$.  The following theorem is a slightly more general version of a result in \cite{D-P1}.

\begin{thm}  \label{structure}
Let $\g$ be a locally reductive Lie algebra.
\begin{enumerate}
\item \label{split} There is a split exact sequence of Lie algebras $$0 \rightarrow
[\g, \g] \rightarrow \g \rightarrow \g/[\g, \g] =: \a \rightarrow 0,$$
i.e. $\g \simeq [\g, \g] \subsetplus \a$, where $\a$ is abelian.

\item If $\g$ is root-reductive, the Lie algebra $[\g, \g]$ \label{simple} is isomorphic to a direct sum of finite-dimensional simple Lie algebras, as well as of copies of $\sl_\infty$, $\so_\infty$, and $\sp_\infty$, with at most countable multiplicities.
\end{enumerate}
\end{thm}

\begin{proof}
\begin{enumerate}
\item Fix an exhaustion $\g = \bigcup_i \g_i$, where $\g_i$ is finite dimensional and reductive.  
Let $\h_i \subseteq \g_i$ be a family of Cartan subalgebras such that $\h_i \cap \g_{i-1} = \h_{i-1}$.  Define $\h := \bigcup_i \h_i$.  Then $\g = [\g,\g] + \h$.  Since $\h$ is abelian, there is an abelian Lie algebra $\a \subseteq \h$ such that $\h = \h \cap [\g,\g] \oplus \a$.  One may check that $\g \simeq [\g,\g] \subsetplus \a$. 
\item For the proof see \cite[Theorem 1]{D-P1}.
\end{enumerate}
\end{proof}

We will apply the following proposition only in the case of root-reductive Lie algebras, but we prove it in more generality for locally reductive Lie algebras.

\begin{prop} \label{extendrep}
Let $\g = \bigcup_i \g_i$ be a locally reductive Lie algebra, with $\g_i$ finite dimensional and reductive, and $\g_i$ reductive in $\g_{i+1}$.  Let $W$ be a $[\g,\g]$-module which is a direct limit of injections of simple $[\g_i,\g_i]$-modules $W_i$.  Then the $[\g,\g]$-module structure on $W$ can be extended to a $\g$-module structure.
\end{prop}

\begin{proof}
For each $i$, choose an abelian Lie algebra $\a_i \subseteq \z(\g_i)$ such that $\g_i = (\g_i \cap [\g,\g]) \oplus \a_i$.  Let $\pi_1 : \g_i \rightarrow \g_i \cap [\g,\g]$ and $\pi_2 : \g_i \rightarrow \a_i$ be the corresponding projections.  Consider the inclusion $\varphi : \g_i \hookrightarrow \g_{i+1}$, and note that $\pi_2 \circ \varphi|_{\a_i} : \a_i \rightarrow \a_{i+1}$ is injective.  Let $\b_{i+1}$ be a vector space complement of $\pi_2 \circ \varphi (\a_i)$ in $\a_{i+1}$.  

Given an action of $\a_i$ on $W_i$ by scalar endomorphisms, we define an action of $\a_{i+1}$ on $W_{i+1}$ by scalars as follows.  
Let $a \in \a_i$ be arbitrary, and let $\alpha$ denote the scalar by which $a$ acts on $W_i$.  Observe that the image of $\pi_1\circ \varphi (a)$ in $\End W$ must commute with the image of $[\g_i, \g_i]$.  Hence the image of $\pi_1 \circ \varphi (a)$ in $\End W$ is the direct sum of a scalar endomorphism of $W_i$ and an endomorphism of a vector space complement of $W_i$.  Let $\beta$ denote the scalar by which $\pi_1 \circ \varphi (a)$ acts on $W_i$.  Define the action of $\pi_2 \circ \varphi (a)$ on $W_{i+1}$ to be by the scalar $\alpha - \beta$.  Since $\pi_2 \circ \varphi|_{\a_i}$ is injective, this procedure defines an action of $\pi_2 \circ \varphi (\a_i)$ on $W_{i+1}$.  Define the action of $\b_{i+1}$ on $W_{i+1}$ to be trivial.

We are now ready to define an action of $\g$ on $W$.  If $x \in \g$ and $w \in W$, then for some $i$ we have $x \in \g_i$ and $w \in W_i$.  Set $x \cdot w := \pi_1(x) \cdot w + \pi_2(x) \cdot w$, where the second action on the right hand side comes from the preceding paragraph.  One must check that this is well-defined, i.e. compatible with the inclusions $\g_i \hookrightarrow \g_{i+1}$ and $W_i \hookrightarrow W_{i+1}$.  Explicitly, the composition $\g_i \rightarrow \g_{i+1} \rightarrow \Hom(W_{i+1}, W) \rightarrow \Hom(W_i,W)$ coincides with the map $\g_i \rightarrow \Hom(W_i,W)$.  Since elements of $\a_i$ act as scalars on $W_i$, $W$ is a well-defined $\g$-module.
\end{proof}

In what follows $\g$ will always denote a root-reductive Lie algebra.  

\subsection{The case where $\so_\infty$ is a direct summand of $[\g,\g]$}\label{sosubsection}

In this subsection $\so (V) = \bigwedge^2 V$ is a direct summand of $[\g,\g]$, and $\t$ is a maximal toral subalgebra of $\g$.  We consider $V$ as a $[\g,\g]$-module in which all direct summands other than $\so(V)$ act trivially.  By Proposition \ref{extendrep}, $V$ may be endowed with the structure of a $\g$-module.  We write $V'$ (and more generally $A'$) for the maximal locally finite $\t$-submodule of $V$ (resp. of a $\t$-module $A$).  The key fact about maximal locally finite submodules is that if $A$ and $B$ are two $\t$-modules then $(A \otimes B)' = A' \otimes B'$.  This is a special case of a result of Dimitrov and Zuckerman appearing in \cite[Appendix Prop. A]{N-P}.  Thus $\so(V)' = (V' \otimes V') \cap \bigwedge^2 V = \bigwedge^2 V'$.

Recall that for any Lie subalgebra $\k \subseteq \g$, $\bigoplus_{\alpha \in \supp V \subseteq \k^*} V^\alpha$ is a generalized weight $\k$-submodule of $V$.  If $\k = \t$, it is easy to verify that $\bigoplus_{\alpha \in \supp V \subseteq \t^*} V^\alpha$ coincides with the maximal locally finite $\t$-submodule $V' \subseteq V$.  Therefore
$$ \so(V)' =  \sum_{\alpha, \beta \in \supp V} V^\alpha \wedge V^\beta,$$
and $\z_{\so(V)}(\t)$ is the $0$-weight space of this weight decomposition of $\so(V)'$.  Thus
$$ \z_{\so(V)}(\t) = {\bigwedge}^2 V^0 \oplus \bigoplus_{\alpha \in I} V^\alpha \wedge V^{-\alpha},$$
where $I \sqcup -I = \supp V \cap - \supp V \setminus \{0\}$.  

\begin{thm}\label{so}
Let $\g$ be a root-reductive Lie algebra for which $\so(V)$ is a direct summand of $[\g, \g]$.  The following statements hold for any maximal toral subalgebra $\t \subseteq \g$:
\begin{enumerate}
\item\label{one} If $\langle V^\alpha , V^\beta \rangle \neq 0$, then $\beta = -\alpha$.  If furthermore $\alpha \neq 0$, then $\dim V^\alpha = \dim V^{-\alpha}  =1$;
\item \label{two} If $\langle \cdot , \cdot \rangle |_{V^\alpha \times V^{-\alpha}} = 0$ then  $V^\alpha \wedge V^{-\alpha} \subseteq \z (\z_{\so(V)} (\t))$; 
\item \label{three} $\langle \cdot , \cdot \rangle|_{V^0 \times V^0}$ has rank 0 or 1, and $(\bigwedge^2 V^0) \cap \t = 0$;
\item \label{four} $\z_{\so(V)} (\t)$ is nilpotent of depth at most $2$.
\end{enumerate}
If $\g = \so(V)$, a subalgebra $\t \subseteq \g$ is maximal toral if and only if (\ref{three}), 
%(\ref{five}), and (\ref{six})
(5) and (6)
%%%%%%%
 hold, where 
\begin{enumerate}
\item[(5)]\label{five} $\supp V= -\supp V$ and  $\dim V^\alpha = 1$ for any $\alpha \in \supp V\setminus \{ 0\}$;
\item[(6)]\label{six} $\t = \bigoplus_{\alpha \in I} V^\alpha \wedge V^{-\alpha}$, where $I \sqcup -I = \supp V \setminus \{ 0 \}$.
\end{enumerate}
\end{thm}

\begin{proof}
\begin{enumerate}
\item Let $v \in V^\alpha$ and $w \in V^\beta$ be such that $\langle v , w \rangle \neq 0$.  If $t \in \t$, then the $\g$-invariance of $\langle \cdot, \cdot \rangle$ yields $\langle t v , w \rangle + \langle v , t w \rangle = 0$.  Thus $(\alpha(t) + \beta (t))\langle v , w \rangle =0$.  Since $\langle v , w \rangle \neq 0$, it follows that $\beta = -\alpha$.  In particular, if $\alpha \neq 0$, then $\langle V^\alpha, V^\alpha \rangle = 0$.

Suppose $\alpha \neq 0$ and without loss of generality that $\langle v,w \rangle = 1$.  Each of $v \otimes w$ and $w \otimes v$ satisfies the equation $x^2 = x$ and therefore is semisimple.  Furthermore, $[ v \otimes w , w \otimes v] = \langle w , w \rangle v \otimes v - \langle v , v \rangle w \otimes w = 0$.  Hence $v \wedge w = v \otimes w - w \otimes v$ is semisimple.  Since $v \wedge w \in \z_\g(\t)$, $t:= v \wedge w$ is contained in $\t$ by maximality.  Then $\alpha(t) v = t \cdot v =  v$ implies $\alpha(t) = 1$. For $v' \in V^\alpha$, we calculate $v' = \alpha(t) v' = t \cdot v' = \langle v' , w \rangle v$.  Therefore $V^\alpha$ is $1$-dimensional.

\item In view of Part (\ref{one}), $V^\alpha$ and $V^{-\alpha}$ are in the radical of the form $\langle \cdot , \cdot \rangle |_{V' \times V'}$ and hence all elements of $\so(V)'$ centralize $V^\alpha \wedge V^{-\alpha}$.

\item  Suppose, for the sake of a contradiction, that the form $\langle \cdot , \cdot \rangle|_{V^0 \times V^0}$ has rank greater than $1$.  Then there exist $v,w \in V^0$ with pairings given by the matrix $\left[
\begin{matrix}
 1 & 0 \\ 0 & 1
 \end{matrix}
\right]$.  As in Part (\ref{one}), we see that $v \wedge w$ is semisimple and hence lies in $\t$.  Since $w$ is in the $0$-weight space of $\t$, we calculate $0 = (v \wedge w) \cdot w = v $, which contradicts the choice of $v$.  Therefore the form $\langle \cdot , \cdot \rangle|_{V^0 \times V^0}$ has rank $0$ or $1$. 

To prove that $(\bigwedge^2 V^0) \cap \t = 0$, it is enough to show that elements of $\bigwedge^2 V^0$ are nilpotent as endomorphisms of $V$.  Let $\sum v_i \wedge w_i \in \bigwedge^2 V^0$ be arbitrary.  If the form $\langle \cdot , \cdot \rangle|_{V^0 \times V^0}$ has rank $0$, one may compute that $(\sum v_i \wedge w_i )^2 = 0$.  If the form $\langle \cdot , \cdot \rangle|_{V^0 \times V^0}$ has rank $1$, one may compute that $(\sum v_i \wedge w_i )^3 = 0$.

\item We have $\z_{\so(V)}(\t) = \bigwedge^2 V^0 \oplus \bigoplus_{\alpha \in I} V^\alpha \wedge V^{-\alpha}$, and since $\bigoplus_{\alpha \in I} V^\alpha \wedge V^{-\alpha}$ is central in $\z_{\so(V)}(\t)$, it is enough to show that $ \bigwedge^2 V^0$ is nilpotent of depth at most $2$.

If $\langle \cdot , \cdot \rangle|_{V^0 \times V^0} = 0$, then $\bigwedge^2 V^0$ is abelian.  Now suppose that $\langle \cdot , \cdot \rangle|_{V^0 \times V^0}$ has rank $1$.  Let $K$ be the radical of the form $\langle \cdot, \cdot \rangle|_{V^0 \times V^0}$.  Let $v$ be any vector in $V^0 \setminus K$, and we may as well assume that $\langle v, v \rangle = 1$.  We have $V^0 = K \oplus \C v$.  To show that $ \bigwedge^2 V^0$ is nilpotent of depth at most $2$, it is enough to show that  $[K \wedge v, K \wedge v] \subseteq \bigwedge^2 K$, since $\bigwedge^2 K$ is central.  For arbitrary $a,b \in K$, we compute
\begin{eqnarray*}
[ a \wedge v , b \wedge v] & = &
[ a \otimes v - v \otimes a , b \otimes v - v \otimes b] \\
& = & (a \otimes v - v \otimes a) \cdot (b \otimes v - v \otimes b) - 
(b \otimes v - v \otimes b) \cdot (a \otimes v - v \otimes a) \\
& = & - a \otimes b + b \otimes a \\
& = & b \wedge a \in {\bigwedge} ^2  K.
\end{eqnarray*}

\item For the first equality, let us compute $\t \cdot V'$ in two different ways.  On the one hand, by definition, $\t \cdot V'  = \bigoplus_{0 \neq \alpha \in \supp V} V^{\alpha}$.  On the other hand, the assumption $\g = \bigwedge^2 V$ yields $\t \subseteq \bigwedge^2 V^0 \oplus \sum_{\alpha \in \supp V \cap -\supp V} V^\alpha \wedge V^{-\alpha}$, and since $\langle V^\alpha, V^\beta \rangle = \{0\}$ for $\beta \neq -\alpha$ ,
we obtain $\t \cdot V' \subseteq \bigoplus_{\alpha \in \supp V \cap -\supp V} V^\alpha$.  
This implies  $\supp V \subseteq \supp V \cap - \supp V$, and hence $\supp V = -\supp V$.

Now consider the fact that $\t \cdot V^\alpha \neq \{0\}$ for $\alpha \neq 0$.  It follows that $\langle V^\alpha, V^{-\alpha} \rangle \neq \{0\}$.  Part (\ref{one}) gives $\dim V^\alpha = 1$  for any $\alpha \in \supp V\setminus \{ 0\}$.

\item By the proof of Part (\ref{one}), we have $V^\alpha \wedge V^{-\alpha} \subseteq \t$ for $0 \neq \alpha \in \supp V$.  This gives the inclusions $$\bigoplus_{\alpha \in I} V^\alpha \wedge V^{-\alpha} \subseteq \t \subseteq  {\bigwedge}^2 V^0 \oplus \bigoplus_{\alpha \in I} V^\alpha \wedge V^{-\alpha}.$$
In Part (\ref{three}) we showed that $(\bigwedge^2 V^{0}) \cap \t = 0$.  Hence  $\t = \bigoplus_{\alpha \in I} V^\alpha \wedge V^{-\alpha}$.

\end{enumerate} 

Conversely, assume that Conditions (\ref{three}),
%(\ref{five}), and (\ref{six})
(5), and (6)
%%%%%%%%%%
 are satisfied.  
Since $\langle V^\alpha , V^\beta \rangle = 0$ for $\alpha \neq -\beta$, the subalgebra $\t$ is a direct sum of $1$-dimensional Lie subalgebras, hence abelian.  If $\alpha \neq 0$, we have $\t \cdot V^\alpha \neq 0$, hence $\langle V^\alpha , V^{-\alpha} \rangle \neq 0$.  Thus $\t$ is spanned by elements of the form $v \wedge w$ where $v \in V^\alpha$, $w \in V^{-\alpha}$, $\alpha \neq 0$, and $\langle v , w \rangle = 1$.  These elements are semisimple as in Part (\ref{one}).  Hence $\t$ is a toral subalgebra of $\g = \bigwedge^2 V$. 
 The centralizer of $\t$ in $\g$ is contained in $\g'$ and coincides with $\t \oplus \bigwedge^2 V^0.$ Condition (\ref{three}) implies that each element in $\bigwedge^2 V^0$ is nilpotent, so $\t$ is maximal toral. 
\end{proof}

\begin{cor}
For any maximal toral subalgebra $\t \subseteq \so(V)$, we have $$ \z_{\so(V)}(\t) = \t \oplus {\bigwedge}^2  V^0, $$ which is a Lie algebra nilpotent of depth at most $2$.
\end{cor}

The result from Theorem \ref{so} (\ref{three}) that $\langle \cdot , \cdot \rangle|_{V^0 \times V^0}$ has rank 0 or 1 enables us to make the following definition.  A maximal toral subalgebra of $\so_\infty$ is \emph{even} if the rank of $\langle \cdot , \cdot \rangle|_{V^0 \times V^0}$ is 0, and \emph{odd} if the rank is 1.  Likewise, a Cartan subalgebra $\h$ is \emph{even} if $\h_{ss}$ is an even maximal toral subalgebra, and  $\h$ is \emph{odd} if $\h_{ss}$ is odd.

There are two types of exhaustions of $\so_\infty$ via root inclusions: as $\bigcup_i \so_{2i}$ and as $\bigcup_i \so_{2i+1}$.  A union of Cartan subalgebras of $\so_{2i}$ yields an even splitting Cartan subalgebra of $\so_\infty$, and a union of Cartan subalgebras of $\so_{2i+1}$ yields an odd splitting Cartan subalgebra of $\so_\infty$.  It is easy to see that any splitting Cartan subalgebra of $\so_\infty$ arises in one of these two ways and is necessarily abelian.

Notice that $\bigwedge^2 V^0$ is the degenerate orthogonal Lie algebra $\so(V^0)$.  For a degenerate symmetric form on a vector space $W$, the Lie algebra $\so(W)$ is nilpotent precisely when the form has rank $0$ or $1$.  Explicitly, the Lie algebra $\so(W)$ is nilpotent and nonabelian exactly when the form has rank $1$ and $\dim W \geq 3$.  Hence any nonabelian Cartan subalgebra of $\so(V)$ must be odd and must have $\dim V^0 \geq 3$.

Here is an example of a nonabelian Cartan subalgebra of $\so(V)$.  Let $\{ e_i \}_{i \in \Z}$ be a basis for $V$ with $\langle e_i , e_j \rangle = \delta_{i,-j}$.  Let $$\t := \bigoplus_{i \geq 3} \C (e_i + e_1) \wedge (e_{-i} + e_{-2}).$$  Observe that $\t \subseteq \so(V)$ is a maximal toral subalgebra, and we have $$\z_{\so(V)}(\t) = \t \oplus {\bigwedge}^2 \Span\{e_0,e_{-2},e_1\}.$$  We see that the Cartan subalgebra $\z_{\so(V)}(\t)$ is not abelian as $[e_1 \wedge e_0, e_{-2} \wedge e_0] = e_{-2} \wedge e_1$.  

It is easy to check that the centralizer of any even maximal toral subalgebra acts locally finitely on itself, whereas the centralizer of an odd maximal toral subalgebra acts locally finitely on itself if and only if $V^0$ is finite dimensional.   

\subsection{The case where $\sp_\infty$ is a direct summand of $[\g,\g]$}

Let $\sp (V)=\Sym^2(V)$ be a direct summand of $[\g,\g]$, and $\t$ a maximal toral subalgebra of $\g$.  The following results are proved in the same way as the analogous statements for $\so(V)$.

We will write $A \& B := \{ a \otimes b + b \otimes a : a \in A, b \in B \} \subseteq A \otimes B \oplus B \otimes A$ for the symmetrizer of vector spaces $A$ and $B$.   Set also $a \& b :=  a \otimes b + b \otimes a  \in A \& B$ for $a \in A$, $b \in B$.

We consider $V$ as a $[\g,\g]$-module in which all direct summands other than $\sp(V)$ act trivally, and by Proposition \ref{extendrep}, $V$ may be endowed with the structure of a $\g$-module.  We have the equality $\sp(V)' = \Sym^2 (V')$ and the weight decomposition $$\sp(V)' = \sum_{\alpha, \beta \in \supp V} V^\alpha \& V^\beta.$$ 
Observe that
$$ \z_{\sp(V)}(\t) = \Sym^2 (V^0) \oplus \bigoplus_{\alpha \in I} V^\alpha \& V^{-\alpha},$$ where $I \sqcup -I = \supp V \cap - \supp V \setminus \{ 0 \}$.

\begin{thm}\label{sp}
Let $\g$ be a root-reductive Lie algebra for which $\sp(V)$ is a direct summand of $[\g, \g]$.  
The following statements hold for any maximal toral subalgebra $\t \subseteq \g$:
\begin{enumerate}
\item If $\langle V^\alpha , V^\beta \rangle \neq 0$, then $\beta = -\alpha$.  If furthermore $\alpha \neq 0$, then $\dim V^\alpha = \dim V^{-\alpha}  =1$;
\item If $\langle \cdot , \cdot \rangle |_{V^\alpha \times V^{-\alpha}} = 0$ then  $V^\alpha \& V^{-\alpha} \subseteq \z (\z_{\sp(V)} (\t))$;  
\item $\langle \cdot , \cdot \rangle |_{V^0 \times V^0}=0$, and $(\Sym^2 (V^0)) \cap \t = 0$;
\item \label{spfour} $\z_{\sp(V)} (\t)$ is abelian.
\end{enumerate}
If $\g=\sp(V)$, a subalgebra $\t \subseteq \g =  \sp(V)$ is maximal toral if and only if  (\ref{three}), 
%(\ref{five}), and (\ref{six})
(5), and (6)
%%%%%%%
 hold, where 
\begin{enumerate}
\item[(5)] $\supp V= -\supp V$ and  $\dim V^\alpha = 1$ for any $\alpha \in \supp V \setminus \{ 0\}$;
\item[(6)] $\t = \bigoplus_{\alpha \in I} V^\alpha \& V^{-\alpha}$, where $I \sqcup -I = \supp V \setminus \{ 0 \}$.
\end{enumerate}
\end{thm}

\begin{cor}
For any maximal toral subalgebra $\t \subseteq \sp(V)$, we have $$\z_{\sp(V)}(\t) = \t \oplus \Sym^2 (V^0),$$  which is an abelian Lie algebra.
\end{cor}

\subsection{The case where $\sl_\infty$ is a direct summand of $[\g,\g]$}

In this subsection we generalize a theorem from \cite{N-P}, which we first recall.
\begin{thm} \cite{N-P} \label{thmc}
A  subalgebra $\t \subseteq \g = \gl(V,V_*)$ is maximal toral if and only if the following conditions are satisfied: 
\begin{enumerate}
\item $(\supp V)\setminus \{ 0\} = -(\supp V_*)\setminus \{ 0\}$ and  $\dim V^\alpha = \dim V_*^{-\alpha} = 1$ for any $\alpha \in \supp V\setminus \{ 0\}$;
\item $\langle \cdot , \cdot \rangle|_{V^0 \times V_*^0} = 0$, and $(V^0 \otimes V_*^0) \cap \t = 0$;
\item $\t = \bigoplus_{0\neq\alpha \in \supp V} V^\alpha \otimes V_*^{-\alpha}$.
\end{enumerate}

Thus, for any maximal toral subalgebra $\t \subseteq \gl(V,V_*)$, we have $$ \z_{\gl(V,V_*)}(\t) = \t \oplus (V^0 \otimes V_*^0), $$
which is an abelian Lie algebra.
\end{thm}

Assume that $\sl(V,V_*)$ is a direct summand of $[\g,\g]$, and let $\t$ be a maximal toral subalgebra of $\g$.  We consider $V$ and $V_*$ as $[\g,\g]$-modules in which all direct summands other than $\sl(V,V_*)$ act trivially, and by Proposition \ref{extendrep}, $V$ and $V_*$ may be endowed with the structure of $\g$-modules.  We have the equality $\sl(V,V_*)' = \sl(V',V_*')$ and the weight decomposition $$\sl(V,V_*)' = \sl(V,V_*) \cap \bigoplus_{\alpha \in \supp V,\beta \in \supp V_*} V^\alpha \otimes V_*^\beta.$$ 
Observe that
$$\z_{\sl(V,V_*)}(\t) = \sl(V,V_*) \cap \bigoplus_{\alpha \in \supp V} V^\alpha \otimes V_*^{-\alpha}.$$ 

\begin{thm}\label{sl}
Let $\g$ be a root-reductive Lie algebra for which  $\sl(V,V_*)$ is a direct summand of $[\g, \g]$.  The following statements hold for any maximal toral subalgebra $\t \subseteq \g$:
\begin{enumerate}
\item\label{slone} If $\langle V^\alpha , V_*^\beta \rangle \neq 0$, then $\beta = -\alpha$.  If furthermore the rank of $\langle \cdot , \cdot \rangle |_{V' \times V_*'}$ is at least $2$, then $\dim V^\alpha = \dim V_*^{-\alpha}  =1$;
\item \label{sltwo} If $\langle \cdot , \cdot \rangle |_{V^\alpha \times V_*^{-\alpha}} = 0$ then  $V^\alpha \otimes V_*^{-\alpha} \subseteq \z (\z_{\sl(V,V_*)} (\t))$; 
\item \label{slthree} If the rank of $\langle \cdot , \cdot \rangle |_{V' \times V_*'}$ is at least $2$, then $\langle \cdot , \cdot \rangle |_{V^0 \times V_*^0} = 0$.  Hence $\langle \cdot , \cdot \rangle |_{V^0 \times V_*^0}$ has rank $0$ or $1$.  Also, $(V^0 \otimes V_*^0) \cap \t = 0$;
\item \label{slfour} $\z_{\sl(V,V_*)} (\t)$ is nilpotent of depth at most $2$.
\end{enumerate}
If $\g = \sl(V,V_*)$, a subalgebra $\t \subseteq \g$ is maximal toral if and only if 
%(\ref{slthreeprime}), (\ref{slfive}), and (\ref{slsix})
(3'), (5), and (6)
%%%%%%%
 hold, where 
\begin{enumerate}
\item[(3')] $\langle \cdot , \cdot \rangle|_{V^0 \times V_*^0} = 0$, and $(V^0 \otimes V_*^0) \cap \t = 0$;
\item[(5)] $(\supp V)\setminus \{ 0\} = -(\supp V_*)\setminus \{ 0\}$ and  $\dim V^\alpha = \dim V_*^{-\alpha} = 1$ for any $\alpha \in \supp V\setminus \{ 0\}$;
\item[(6)] $\t =  \sl(V,V_*) \cap \bigoplus_{0\neq\alpha \in \supp V} V^\alpha \otimes V_*^{-\alpha}$.
\end{enumerate}
\end{thm}

\begin{proof}
\begin{enumerate}
\item Suppose $\langle V^\alpha , V_*^{-\alpha} \rangle \neq 0$ and that the rank of $\langle \cdot , \cdot \rangle |_{V' \times V_*'}$ is at least $2$.  Then let $v \in V^\alpha$ and $w \in V_*^{-\alpha}$ be such that $\langle v , w \rangle =1$.  Since the rank of  $\langle \cdot , \cdot \rangle |_{V' \times V_*'}$ is at least $2$, there exist  $v' \in V^\beta$ and $w' \in V_*^{-\beta}$ such that $\langle v' , w' \rangle =1$ and  $\langle v , w' \rangle = \langle v' , w \rangle =0$.  We have $t := v \otimes w - v' \otimes w' \in \t$, by an argument similar to the $\so_\infty$ case.  Then $t \cdot v = v$ implies $\alpha(t) = 1$.  

First assume $\alpha \neq \beta$.  Let $v'' \in V^\alpha$, and compute $v'' = t \cdot v'' = (v \otimes w - v' \otimes w') \cdot v'' = \langle v'' , w \rangle v$.  This implies that $V^\alpha$ is $1$-dimensional.  A similar argument shows that $V_*^{-\alpha}$ is also $1$-dimensional.

Second assume, for the sake of a contradiction, that $\alpha = \beta$.  Then compute $t \cdot v' = -v'$, which implies $\alpha(t) = -1$, contradicting  $\alpha(t) = 1$.

\item The argument from the $\so_\infty$ case works; see Theorem \ref{so}.
\item If the rank of $\langle \cdot , \cdot \rangle |_{V' \times V_*'}$ is at least $2$, and if $\langle V^\alpha, V_*^{-\alpha} \rangle \neq 0$, then as in Part (\ref{slone}) there exists $t \in \t$ with $\alpha (t) =1$.  Hence   $\langle V^0, V_*^0 \rangle = 0$.  The proof of the last statement is similar to that in Theorem \ref{so}.

\item We have $\z_{\sl(V,V_*)}(\t) = \sl(V,V_*) \cap \bigoplus_{\alpha} V^\alpha \otimes V_*^{-\alpha}$.  If the rank of $\langle \cdot , \cdot \rangle |_{V' \times V_*'}$ is not $1$, then one may check that $\z_{\sl(V,V_*)}(\t)$ is abelian.  Suppose the rank of $\langle \cdot , \cdot \rangle |_{V' \times V_*'}$ is $1$.  Let $\beta$ be the weight for which $\langle \cdot , \cdot \rangle |_{V^\beta \times V_*^{-\beta}}$ has rank $1$.  Then $\bigoplus_{\alpha \neq \beta} V^\alpha \otimes V_*^{-\alpha}$ is central in $\z_{\sl(V,V_*)}(\t)$.  Therefore it suffices to show that $\sl(V^\beta, V_*^{-\beta})$ is nilpotent of depth at most $2$.  It is an easy computation to check that if either $V^\beta$ or $V_*^{-\beta}$ has dimension $1$, then $\sl(V^\beta, V_*^{-\beta})$ is abelian, and that otherwise it is nilpotent of depth $2$.
\end{enumerate}

If the rank of $\langle \cdot , \cdot \rangle |_{V' \times V_*'}$ is $0$ or $1$, then $\t \cap \sl(V,V_*) = 0$.  Thus for any maximal toral subalgebra $\t \subseteq \sl(V,V_*)$, we have that $\langle \cdot , \cdot \rangle |_{V^0 \times V_*^0}$ has rank $0$, which yields (3').  One can modify the proofs from Theorem \ref{so} to obtain (5) and (6), as well as the reverse implication.
%%%%%%%%%%%%%%BAD REFERENCING
\end{proof}

\begin{cor}
For any maximal toral subalgebra $\t \subseteq \sl(V,V_*)$, we have $$ \z_{\sl(V,V_*)}(\t) = \t \oplus (V^0 \otimes V_*^0),$$ which is an abelian Lie algebra.
\end{cor}

Note that the degenerate Lie algebra $\sl(W,W_*)$ is nilpotent if and only if the form $\langle \cdot , \cdot \rangle : W \times W_* \rightarrow \C$ has rank $0$ or $1$.  In the rank $1$ case, the algebra is again nilpotent of depth at most $2$.  This phenomenon does not occur for $\gl_\infty$ or $\sp_\infty$, since in these cases any nilpotent degenerate algebra must be fully degenerate. 

\subsection{Unified description of Cartan subalgebras of $\gl_\infty$, $\sl_\infty$, $\so_\infty$, and $\sp_\infty$}

\begin{defn}
\begin{enumerate}
\item A \emph{dual system for $\gl(V, V_*)$} or \emph{for $\sl(V,V_*)$} is a set of of $1$-dimensional vector subspaces $L_i \subseteq V$ and $L^i \subseteq V_*$ for $i \in I$ such that $\langle L_i, L^j \rangle = \delta_{ij} \C$.  Here $I$ is a finite set or $Z_{>0}$.

\item A \emph{self-dual system for $\so(V)$} or \emph{for $\sp(V)$} is a set of $1$-dimensional vector subspaces $L_i \subseteq V$ for $i \in I$ such that $\langle L_i, L_j \rangle = \delta_{i, -j} \C$. Here $I$ is a finite subset of $\Z_{\neq 0}$ with $-I = I$, or $I = \Z_{\neq 0}$.
\end{enumerate}
\end{defn}

The following is a corollary to Theorems  \ref{so}, \ref{sp}, \ref{thmc}, and \ref{sl}.

\begin{cor}\label{bijection}
Cartan subalgebras of $\gl_\infty$ (resp. of $\sl_\infty$, $\so_\infty$, or $\sp_\infty$) are in one-to-one correspondence with maximal (self-)dual systems for $\gl_\infty$ (resp. for $\sl_\infty$, $\so_\infty$, or $\sp_\infty$).
\end{cor}

\begin{proof} 
A (self-)dual system yields a toral subalgebra by the formulas in Figure \ref{fig1}.  In the case of $\gl_\infty$, since $\langle L_i , L^j \rangle = 0$ for $i \neq j$, the subalgebra $\t$ is a direct sum of $1$-dimensional Lie subalgebras, hence abelian. Since $\langle L_i , L^i \rangle \neq 0$, $\t$ is spanned by elements of the form $v \otimes w$, where $v \in L_i$, $w \in L^i$, and $\langle v , w \rangle = 1$.  These elements are semisimple since they satisfy the equation $x^2 = x$.  Hence $\t$ is a toral subalgebra.  A similar argument shows that self-dual systems also yield toral subalgebras.
 
Consider the map which sends a (self-)dual system to the toral subalgebra according to the formulas in Figure \ref{fig1}.  This map is injective (except in the case of $\sl_\infty$, where one must not allow dual systems to have $|I| =1$).  To see injectivity in the case of $\gl_\infty$, suppose $\{L_i,L^i \}$ and $\{ M_j , M^j \}$ are dual systems with $\bigoplus_i L_i \otimes L^i = \t = \bigoplus_j M_j \otimes M^j$.  Then $L_i = \t \cdot L_i = (\bigoplus_j M_j \otimes M^j) \cdot L_i$.  Hence for some $j$ we have $M_j = L_i$.  This argument can be adapted to work for self-dual systems.
Moreover, this map preserves containment, and its image includes all maximal toral subalgebras by Theorems \ref{so}, \ref{sp}, \ref{thmc}, and \ref{sl}.  Therefore maximal toral subalgebras correspond to maximal (self-)dual systems.
\end{proof}

Explicitly, the maximal toral subalgebra $\t$ and the Cartan subalgebra $\h$ in $\g$ associated to a maximal (self-)dual system for $\g$ are given in Figure \ref{fig1}.  Note that a maximal toral subalgebra $\t \subseteq \g$ is splitting if and only if $V=V'$ and $V_* = V_*'$.  In that case necessarily $(\bigoplus_i L_i )^\bot = 0$ and $(\bigoplus_i L^i )^\bot = 0$ for  $\gl_\infty$ or $\sl_\infty$, $(\bigoplus_i L_i )^\bot = 0$ for $\sp_\infty$, and for $\so_\infty$ $(\bigoplus_i L_i )^\bot$ is $0$ or $1$-dimensional.  In fact, these conditions are equivalent to being splitting.

\begin{figure}
\caption{Maximal toral subalgebra $\t$ and Cartan subalgebra $\h$ associated to a maximal (self-)dual system}
 \label{fig1}
$$\begin{array}{c|c|c}
\g & \t & \h \\ \hline
\gl_\infty & \bigoplus_i L_i \otimes L^i & \t \oplus (( \bigoplus_i L^i )^\bot \otimes ( \bigoplus_j L_j )^\bot) \\ \hline
\sl_\infty &  \sl_\infty \cap (\bigoplus_i L_i \otimes L^i) & \t \oplus (( \bigoplus_i L^i )^\bot \otimes ( \bigoplus_j L_j )^\bot) \\ \hline
\so_\infty & \bigoplus_{i>0} L_i \wedge L_{-i} & \t \oplus \bigwedge^2 \big( (\bigoplus_i L_i)^\bot \big) \\ \hline
\sp_\infty & \bigoplus_{i>0} L_i \& L_{-i} & \t \oplus \Sym^2 \big( (\bigoplus_i L_i)^\bot \big)
\end{array}$$
\end{figure}

\subsection{The case of a general $\g$.}\label{nilp}

The following theorem is our main result on Cartan subalgebras of general root-reductive Lie algebras and strengthens Theorem \ref{main} in this case.

\begin{thm}\label{depth2}
Let $\g$ be a root-reductive Lie algebra and $\h \subseteq \g$ a Cartan subalgebra.  Then $\h$ is nilpotent of depth at most $2$, and $\h = \g^0(\h)$.
\end{thm}

\begin{proof}
Let $\t := \h_{ss}$, so $\h = \z_\g(\t)$.  We first show that $\z_{[\g,\g]} (\t)$ is nilpotent of depth at most $2$.  Recalling the decomposition of $[\g,\g]$ in Theorem \ref{structure} (\ref{simple}), let $\g_{fd}$ be the direct sum of all finite-dimensional simple direct summands of $[\g,\g]$, and $\g_{id}$ the direct sum of all infinite-dimensional simple direct summands of $[\g,\g]$.  As there are no nontrivial extensions of an abelian Lie algebra by a finite-dimensional reductive Lie algebra, we get the decomposition $$\g = \g_{fd} \oplus ( \g_{id} \subsetplus \a),$$
where $\a$ is abelian.  Then $\t = \t_1 \oplus \t_2$ where $\t_1$ is a maximal toral subalgebra of $\g_{fd}$ and $\t_2$ is a maximal toral subalgebra of $\g_{id} \subsetplus \a$.  So $\z_{[\g,\g]}(\t) = \z_{\g_{fd}} (\t_1) \oplus \z_{\g_{id}} (\t_2)$.
As $\t_1$ is self-centralizing in $\g_{fd}$, it remains to show that $\z_{\g_{id}} (\t_2)$ is nilpotent of depth at most $2$.  Let $\g_{id} = \oplus_j \s_j$ be the decomposition of $\g_{id}$ into simple direct summands.  Since $[\s_i, \g] \subseteq \s_i$, we have
 $$\z_{(\oplus_j \s_j)} (\t_2) = \oplus_j \z_{\s_j} (\t_2).$$
Theorems \ref{so} (\ref{four}), \ref{sp} (\ref{spfour}), and \ref{sl} (\ref{slfour}) imply that  $\z_{\g_{id}} (\t_2)$ is nilpotent of depth at most $2$.  Thus $\z_{[\g, \g]} (\t)$ is a Lie algebra nilpotent of depth at most $2$.

Recall from Theorem \ref{main} that $\z_\g (\t)$ is splittable.
If $h = h_{ss} + h_{nil} \in \z_\g (\t)$, then $h_{ss} , h_{nil} \in \z_\g(\t)$.  By definition $h_{ss} \in \t$, and $h_{nil} \in [\g, \g]$ as $h_{nil}$ is nilpotent.  Thus $\z_\g (\t) = \t + \z_{[\g, \g]} (\t) $.  Since $\t$ is central in $\z_\g (\t)$, it follows that $\z_\g (\t)$ is nilpotent of depth no greater than $2$.

Finally, the nilpotence of $\h$ implies $\h \subseteq \g^0(\h)$.  Therefore the condition $\h = \overline{\g^0(\h)}$ together with the inclusions $\h \subseteq \g^0(\h) \subseteq \overline{\g^0(\h)}$ yields $\h = \g^0(\h)$.
\end{proof}

 In the case of $\g = \gl_\infty$, any maximal toral subalgebra of $\g$ surjects onto
$\g/[\g,\g]$ \cite{N-P}.  In general, the map $\t \rightarrow
\g/[\g,\g]$ does not have to be surjective.  For example, consider the Lie algebra $\g$ defined as the direct limit of the inclusions
\begin{eqnarray*}
\gl_{2n} & \hookrightarrow & \gl_{2n+2} \\
A  & \mapsto &
\left[
\begin{matrix} 
\frac{\tr(A)}{n} & & &\\
& A & &\\
& & & 0
\end{matrix}
\right] .
\end{eqnarray*}
Clearly we have the exact sequence 
$$0 \rightarrow \sl_\infty = [\g,\g] \rightarrow \g \rightarrow \C \rightarrow 0.$$
Let $B_{2n} \in \gl_{2n}$ denote the matrix 
$$ 
\left[
\begin{matrix}
I_n & 0_n \\
0_n & 0_n 
\end{matrix}
\right],$$
and let $B$ denote the element of $\g$ defined by the sequence $(B_{2n})$.  
Then as a vector space $\g$ is isomorphic to $\sl_\infty \oplus \C B$.  We will exhibit a Cartan subalgebra of $\g$ whose image in $\C$ is trivial.

For $n \in \Z_{>0}$, let $C_n \in \sl_{2n}$ be the matrix whose only nonzero entry is a $1$ in the upper right corner, and let $C_{-n} \in \sl_{2n}$ be the matrix whose only nonzero entry is a $1$ in the lower left corner.  Then $\t := \Span_{n > 0} \{ C_n + C_{-n} \}$ is a toral subalgebra of $\g$.  We will show that  $\z_\g(\t) \subseteq \sl_\infty$.  A general element of $\g$ lies in $\gl_{2n}$ for some $n > 0$, and hence it may be expressed as $M + aB$ for some $M \in \sl_{2n}$ and $a \in \C$.  The computation $[M +aB , C_{n+1} + C_{-(n+1)}] = a(C_{n+1} - C_{-(n+1)})$ implies that if $M + aB$ centralizes $\t$, then $a=0$.  Let $\h$ be a Cartan subalgebra which is the centralizer of any maximal toral subalgebra containing $\t$.  The containments $\h \subseteq \z_\g(\t) \subseteq \sl_\infty$ imply that $\h$ maps trivially in $\C$.

Theorem \ref{depth2} also leaves open the question whether the intersection $\t \cap [\g,\g]$ is in general a maximal toral subalgebra of $[\g,\g]$.  Clearly in the case of $\gl_\infty$ the intersection of a maximal toral subalgebra with $\sl_\infty$ is maximal.  As the following example shows, the intersection is not in general maximal, and it can even be trivial.

Let $\{ e_i : i \in \Z_{>0} \}$ be a basis of $V$, and let $V_j := \Span \{e_1, \ldots e_j \}$.  Define constants $c_{ij}$ for $i,j \in \Z_{>0}$ by
$$c_{ij} := 
\begin{cases}
1 & j \equiv i (\mathrm{mod} \; 2^{i-1}) \\
0 & \text{otherwise.}
\end{cases}
$$
Let $f_j := \sum_{i \in \Z_{>0}} c_{ij} e_i,$ and notice that $f_j$ is a well-defined vector in $V_j$ as $c_{ij} = 0$ for $i > j$.  Observe also that $\{ f_1 , \ldots f_j \}$ is a basis of $V_j$.  We construct a root-reductive Lie algebra $\g$ with $[\g,\g] = \sl(V,V_*)$ and $\g / [\g,\g]$ countable dimensional.  We consider the maximal toral subalgebra in $\g$ consisting of all elements of $\g$ which have all the $f_j$'s as eigenvectors.  No nontrivial element of $\sl(V,V_*)$ has this property.  

By construction, $f_j - e_j = f_i$ for some unique $i < j$.  Let $p : \Z_{\geq 2} \rightarrow \Z_{\geq 1}$ be defined by $f_j - e_j = f_{p(j)}$.  In addition, define inductively positive integers $d_k^l$ for $k, l \in \Z_{>0}$ by
$$d_k^l := 
\begin{cases}
k & 1 \leq k \leq l \\
d_{p(k)}^l & k > l.
\end{cases}
$$
Let $\g_i := \sl(V_i,V_i^*) \oplus \C^i$, and let $\g$ be the root-reductive Lie algebra defined as the direct limit of the inclusions
\begin{eqnarray*}
\g_i & \hookrightarrow & \g_{i+1} \\
\left( A,0 \right) & \mapsto & \left( \left[
\begin{matrix} 
A & \\
 & 0
\end{matrix}
\right] ,0 \right)\\
\left( 0, \frac{\sum_{k=1}^i d_k^l}{i}x_l \right) & \mapsto &  \left( \left[
\begin{matrix} 
-\frac{\gamma_i^l}{i}I_i & \\
 & \gamma_i^l
\end{matrix}
\right] ,
 \frac{\sum_{k=1}^{i+1} d_k^l}{i+1}x_l \right),
\end{eqnarray*}
where $\{x_l : 1 \leq l \leq i \}$ is a basis of $\C^i$, $I_i$ is the identity matrix in $\gl(V_i,V_i^*)$, and $\gamma_i^l := (i \cdot d_{i+1}^l - \sum_{k=1}^i d_k^l)/(i+1)$.
Clearly $[\g,\g] = \sl(V,V_*)$, where $V_* := \bigcup_i V_i^*$.

For $l \in \Z_{>0}$, define the element $t_l \in \g_l$ by
 $$t_l := \left( (c_{ij})_{1 \leq i,j \leq l} 
\left( \left[
\begin{matrix} 
1 & & \\
&  \ddots &\\
& & l
\end{matrix}
\right]
 - \frac{l+1}{2}I_l \right)
  (c_{ij})_{1 \leq i,j \leq l} ^{-1},
  \frac{l+1}{2}x_l \right).$$
Since $t_l$ is the sum of a central element and the conjugate of a diagonal element of $\sl(V_l,V_l^*)$, $t_l$ is semisimple.  Because $t_l$ and $t_m$ have the same set of eigenvectors, they commute.  Explicitly, if $m \geq l$, the image of $t_l$ in $\g_m$ is
$$\left( (c_{ij})_{1 \leq i,j \leq m} 
\left( \left[
\begin{matrix} 
d_1^l & & \\
&  \ddots &\\
& & d_m^l
\end{matrix}
\right]
 - \frac{\sum_{k=1}^m d_k^l}{m} I_m \right)
  (c_{ij})_{1 \leq i,j \leq m} ^{-1},
 \frac{\sum_{k=1}^m d_k^l}{m}  x_l \right),$$
 i.e. $f_1, \ldots, f_m$ are eigenvectors of the image of $t_l$ in $\g_m$.
Let $\t$ be the toral subalgebra of $\g$ generated by the elements $t_l$.  We will show that no nontrivial element of $\sl(V,V_*)$ centralizes $\t$.

Suppose $C \in \sl(V,V_*)$ centralizes $\t$.  Since $C \in \sl(V_j,V_j^*)$ for some $j$, we may consider $C$ as $
\left( \left[
\begin{matrix} 
C & 0 \\
0 & 0
\end{matrix}
\right],0 \right) \in \g_l$
for any $l > j$. 
In the same block notation, write $t_l = \left( \left[
\begin{matrix} 
* & S_l \\
0 & *
\end{matrix}
\right],* \right)$.
Then $$0 = \left[
\left( \left[
\begin{matrix} 
C & 0 \\
0 & 0
\end{matrix}
\right], 0 \right) ,
\left( \left[
\begin{matrix} 
* & S_l \\
0 & *
\end{matrix}
\right],* \right)
\right] = 
\left( \left[
\begin{matrix} 
* & CS_l \\
0 & 0
\end{matrix}
\right],0 \right) .
$$
The $k$-th column of $S_l$ is $(j+k-p(j+k))\sum_{i=1}^j c_{i,j+k}e_i$.   Let $l$ be sufficiently large that the vectors $\{ \sum_{i=1}^j c_{i,j+k} e_i : 1 \leq k \leq l-j \}$ span $V_j$. 
Since the columns of $S_l$ span $V_j$, it follows that $C = 0$.  Hence no maximal toral subalgebra containing $\t$ intersects $\sl(V,V_*)$ nontrivially.

\section{Conjugacy of Cartan subalgebras of simple root-reductive Lie algebras}\label{sectionc}

We define a toral subalgebra $\t$ of  $\g$, where $\g$ is one of $\gl_\infty$, $\sl_\infty$, $\so_\infty$, or $\sp_\infty$, to be \emph{submaximal} if $\t$ is associated to a not necessarily maximal (self-)dual system via the formulas in Figure \ref{fig1}.  Any nonzero submaximal toral subalgebra is associated to a unique (self-)dual system.  Clearly not every toral subalgebra of $\g$ is submaximal, but any maximal toral subalgebra of $\g$ is submaximal, as seen in Corollary \ref{bijection}.  

Let $\GL(V, V_*)$ be the subgroup of $\GL(V)$ preserving $V_* \subseteq V^*$.  In terms of dual bases for $V$ and $V_*$, elements of  $\GL(V, V_*)$ are matrices with finitely many nonzero entries in each row and each column.  Observe that $\gl(V,V_*) = V \otimes V_*$ is a representation of $\GL(V,V_*)$, as is $\sl(V,V_*)$.  Similary, each of $\so(V)$ and $\sp(V)$ is a representation of the corresponding subgroup of $\GL(V,V)$ preserving the bilinear form.  These subgroups are denoted $\SO(V)$ and $\SP(V)$.  In what follows we describe submaximal toral subalgebras of $\gl(V, V_*)$ (resp. of $\so(V)$ or $\sp(V)$) up to conjugation by the group $\GL(V, V_*)$ (resp. by $\SO(V)$ or $\SP(V)$).

\begin{prop}
Let $\g$ be one of $\gl_\infty$, $\sl_\infty$, $\so_\infty$, or $\sp_\infty$, and $G$ its corresponding group defined above.  Two finite-dimensional submaximal toral subalgebras of $\g$ are conjugate by an element of $G$ if and only if they have the same dimension.
\end{prop}

\begin{proof}
Let $\t_1$ and $\t_2$ be finite-dimensional submaximal toral subalgebras of $\g$ of the same dimension.  Then $\t_1,\t_2 \subseteq \g_i$ for some $i$, where $\g = \bigcup_i \g_i$ is an exhaustion by finite-dimensional reductive Lie algebras under root inclusions $\g_i \subseteq \g_{i+1}$.  It is clear that $\t_1$ and $\t_2$, being submaximal, are conjugate in $\g_i$ by an element of the classical algebraic group $G_i$ associated to $\g_i$. Hence they are conjugate in $\g$, since there is an obvious injective homomorphism of $G_i$ into $G$.
\end{proof}

In what follows, we assume that all (self-)dual systems we consider are infinite.  Their corresponding submaximal toral subalgebras will be infinite dimensional.  Clearly all maximal (self-)dual systems are infinite.

\begin{defn}
\begin{enumerate}
\item A dual system $\{L_i, L^i \}$ for $\gl(V,V_*)$ (resp. for $\sl(V,V_*)$) with $\langle \cdot , \cdot \rangle : V \times V_* \rightarrow \C$ and a dual system $\{M_i, M^i \}$ for $\gl(W,W_*)$ (resp. for $\sl(W,W_*)$) with $\langle \cdot , \cdot \rangle' : W \times W_* \rightarrow \C$ are \emph{equivalent} if there exist isomorphisms $\varphi: V \rightarrow W$ and $\varphi: V_* \rightarrow W_*$ and a bijection $\sigma: \Z_{>0} \rightarrow \Z_{>0}$ such that $\langle \cdot, \cdot \rangle = \varphi^* \langle \cdot, \cdot \rangle'$, $\varphi(L_i) =  M_{\sigma(i)}$, and $\varphi(L^i) = M^{\sigma(i)}$.

\item A self-dual system $\{L_i \}$ for $\so(V)$ (resp. for $\sp(V)$) with $\langle \cdot , \cdot \rangle : V \times V \rightarrow \C$ and a self-dual system $\{M_i \}$ for $\so(W)$ (resp. for $\sp(W)$) with $\langle \cdot , \cdot \rangle' : W \times W \rightarrow \C$ are \emph{equivalent} if there exist an isomorphism $\varphi: V \rightarrow W$ and a bijection $\sigma: \Z_{\neq 0} \rightarrow \Z_{\neq 0}$ such that $\langle \cdot, \cdot \rangle = \varphi^* \langle \cdot, \cdot \rangle'$ and $\varphi(L_i) =  M_{\sigma(i)}$.
\end{enumerate}
\end{defn}
It is clear that equivalent (self-)dual systems for the same algebra are precisely those which are conjugate.

\begin{lemma}\label{lemma3}
Any dual system $\{M_i, M^i \}$ for $\gl(W,W_*)$ is equivalent to a dual system for $\gl(V,V_*)$. Similarly, any self-dual system $\{M_i \}$ for $\so(W)$ (resp. for $\sp(W)$) is equivalent to a self-dual system for $\so(V)$ (resp. for $\sp(V)$).
\end{lemma}

\begin{proof}
The first statement follows directly from Mackey's result \cite[p. 171]{Mackey}, which implies the existence of isomorphisms $W \rightarrow V$ and $W_* \rightarrow V_*$ which send $\langle \cdot , \cdot \rangle'$ to $\langle \cdot, \cdot \rangle$.  Under these isomorphisms, the dual system $\{M_i, M^i \}$ maps to an equivalent dual system.  Similarly, any nondegenerate (symmetric or antisymmetric) form on a countable dimensional vector space can be diagonalized, hence there exists an isomorphism $W \rightarrow V$ which sends $\langle \cdot , \cdot \rangle'$ to $\langle \cdot, \cdot \rangle$.  This isomorphism sends $\{M_i \}$ to an equivalent self-dual system.
\end{proof}

Given a dual system $\{L_i,L^i \}$ for $\gl(V,V_*)$, consider the following construction.  Let $X$ be a vector space complement in $V$ of $\bigoplus_i L_i$, and $Y$ a vector space complement in $V_*$ of $ \bigoplus_i L^i$.  Choose for each $i \in \Z_{>0}$ a generator $v_i \in L_i$, and let $v^i \in L^i$ be the vector such that $\langle v_i , v^i \rangle = 1$.  Forgetting the dual system, we can consider the bilinear form $\langle \cdot , \cdot \rangle$ restricted to $X \times Y$, together with the linear functionals $\lambda_i := \langle \cdot, v^i\rangle : X \rightarrow \C$ and $\mu_i := \langle v_i, \cdot \rangle : Y \rightarrow \C$.  This information encodes the entire structure of the dual system, together with $V$, $V_*$, and $\langle \cdot , \cdot \rangle$, up to isomorphism.  Hence we suggest the following as a useful concept.

\begin{defn}
\begin{enumerate}
\item A \emph{complement datum for $\gl_\infty$ or $\sl_\infty$} is a pair of at most countable dimensional vector spaces $X$ and $Y$, a bilinear map $\omega: X \times Y \rightarrow \C$, a sequence of linear functionals $(\lambda_i)$ on $X$, and a sequence of linear functionals $(\mu_i)$ on  $Y$, for $i \in \Z_{>0}$.  

\item A \emph{complement datum for $\so_\infty$} (resp. for $\sp_\infty$)  is an at most countable dimensional vector space $X$ together with a symmetric (resp. antisymmetric) bilinear form $\omega$ on $X$ and a sequence of linear functionals $(\lambda_i)$ on $X$, for $i \in \Z_{\neq 0}$.
\end{enumerate}
\end{defn}

From a dual system for $\gl(V,V_*)$ and choices of $v_i$, $v^i$, $X$ and $Y$ as described above, we can produce a complement datum for $\gl_\infty$.  Such a dual system and complement datum are said to be \emph{compatible}.   Similarly, a dual system for $\sl(V,V_*)$ and a complement datum for $\sl_\infty$ obtained from it by the same procedure are said to be \emph{compatible}.  We will think of compatibility as a relation.

The notion of compatibility between a self-dual system and a complement datum for $\so_\infty$ or $\sp_\infty$ is similar.  Fix a self-dual system $\{L_i \}$ for $\so(V)$ or $\sp(V)$, and we construct a complement datum as follows.  Let $X$ be a vector space complement of $\bigoplus_i L_i$ in $V$, and choose a nonzero $v_i \in L_i$ for $i > 0$, and let $v_{-i} \in L_{-i}$ be determined by $\langle v_i , v_{-i} \rangle = 1$.  Consider the restriction of $\langle \cdot , \cdot \rangle$ to $X \times X$ together with the linear functionals on $X$ given by $\lambda_i := \langle \cdot , v_i \rangle : X \rightarrow \C$.  Again, the self-dual system and any complement datum obtained from it in this way are said to be \emph{compatible}.

We will say that two sequences $a_i$ and $b_i$ are \emph{almost equal}, and write $(a_i) \approx (b_i )$, if $a_i = b_i$ for all but finitely many $i$.  Let $X_0$ be the set of $x \in X$ such that $( \lambda_i(x) ) \approx 0$.  In the case of $\gl_\infty$ or $\sl_\infty$, let $Y_0$ be the set of $y \in Y$ such that $( \mu_i(y) ) \approx 0$.  

For a fixed (self-)dual system, one might choose a compatible complement datum in such a way that $X$ would actually contain the set of $x \in V$ such that $\lambda_i (x) = 0$ for all $i$.  In this case, $X_0$ would be this orthogonal complement.  The examples we give will all have this property.  We do not require this property in general because it is not necessary and complicates the statements and proofs.

Given a complement datum $(X, Y, \omega, (\lambda_i), (\mu_i))$ for $\gl_\infty$ or $\sl_\infty$, consider the formal sum $$\tilde{\omega} := \omega - \sum_{i \in \Z_{>0}} \lambda_i \otimes \mu_i,$$ where in what follows we consider $\omega$ as a linear map $X \otimes Y \rightarrow \C$.  Note that $\tilde{\omega}$ is a well-defined linear map when restricted to $X_0 \otimes Y$ or $X \otimes Y_0$. 

Given a complement datum $(X, \omega, (\lambda_i))$ for $\so_\infty$, we define analogously $$\tilde{\omega} := \omega - \sum_{i \in \Z_{> 0}} \lambda_i \& \lambda_{-i}.$$  In the case of $\sp_\infty$, we define $$\tilde{\omega} := \omega - \sum_{i \in \Z_{>0}} \lambda_i \wedge \lambda_{-i}.$$  Observe that $\tilde{\omega}$ is a linear map when restricted to $X_0 \otimes X$ or $X \otimes X_0$.  Moreover, we see that $\tilde{\omega}$ retains from $\omega$ the property of symmetry or antisymmetry.

\begin{defn}
\begin{enumerate}
\item
A complement datum  $(X, Y, \omega, (\lambda_i), (\mu_i))$ for $\gl_\infty$ or $\sl_\infty$ is \emph{nondegenerate} if for any nonzero $x_0 \in X_0$ there exists $y \in Y$ such that $\tilde{\omega}(x_0, y) \neq 0$, and for any nonzero $y_0 \in Y_0$ there exists $x \in X$ such that $\tilde{\omega}(x,y_0) \neq 0$.  Moreover  $(X, Y, \omega, (\lambda_i), (\mu_i))$  is \emph{maximal} if $\tilde{\omega}$ restricted to $X_0 \otimes Y_0$ is trivial.  

\item
A complement datum $(X, \omega, (\lambda_i))$ for $\so_\infty$ is \emph{nondegenerate} if for any nonzero $x_0 \in X_0$ there exists $x \in X$ such that $\tilde{\omega}(x_0, x) \neq 0$.  Moreover $(X, \omega, (\lambda_i))$ is \emph{maximal} if $\tilde{\omega}$ restricted to $X_0 \otimes X_0$ has rank $0$ or $1$.

\item
A complement datum $(X, \omega, (\lambda_i))$ for $\sp_\infty$ is \emph{nondegenerate} if for any nonzero $x_0 \in X_0$ there exists $x \in X$ such that $\tilde{\omega}(x_0, x) \neq 0$.  Moreover $(X, \omega, (\lambda_i))$ is \emph{maximal} if $\tilde{\omega}$ restricted to $X_0 \otimes X_0$ is trivial.
\end{enumerate}
\end{defn}

\begin{prop}\label{compatibility}
Any complement datum which is compatible with a (self-)dual system is nondegenerate.  For any compatible pair, the complement datum is maximal if and only if the (self-)dual system is maximal.
\end{prop}

\begin{proof}
We prove the proposition in the case of $\gl_\infty$, and the other cases are similar.

Suppose $(X, Y, \omega, (\lambda_i), (\mu_i))$ is a complement datum compatible with a dual system for $\gl(V,V_*)$.  Take a nonzero vector $x_0 \in X_0$.  Since $( \lambda_i(x_0) ) \approx 0$, we see that $x_0-\sum \lambda_i(x_0) v_i$ is a well-defined nonzero vector in $V$.  Hence, by the nondegeneracy of $\langle \cdot, \cdot \rangle$, there exists $w \in V_*$ such that $\langle x_0-\sum \lambda_i(x_0) v_i, w \rangle \neq 0$.  Let $w = y + \sum b_j v^j$, with $y \in Y$.  We compute:
\begin{align*}
0 \neq 
\langle x_0-\sum \lambda_i(x_0) v_i, y + \sum b_j v^j \rangle
&= \omega(x_0,y) - \sum \lambda_i(x_0) \mu_i(y) \\
& = \tilde{\omega}(x_0,y).
\end{align*}
Similarly, if $y_0 \in Y_0$ the analogous calculation shows that there exists $x \in X$ such that $\tilde{\omega}(x,y_0) \neq 0$.

We turn to maximality.  Suppose that the dual system $\{\C v_i, \C v^i \}$ for $\gl(V,V_*)$ is not maximal.  Then there exist vectors $v_0 \in V$ and $v^0 \in V_*$ such that $\langle v_i, v^j \rangle = \delta_{ij}$.  Let $v_0 = x + \sum a_j v_j$ and $v^0 = y + \sum b_k v^k$ for $x \in X$ and $y \in Y$.  Note $( a_j ) \approx 0$ and $( b_j ) \approx 0$.  We calculate $\lambda_i(x) = \langle x, v^i \rangle = -a_i$, and $\mu_i(y) = -b_i$.  Hence $( \lambda_i(x) ) \approx 0$ and $( \mu_i(y) ) \approx 0$, so $x \in X_0$ and $y \in Y_0$.  Now, 
\begin{align*}
1 &= \langle v_0, v^0 \rangle = \langle x, y \rangle + \sum b_k \lambda_k(x) + \sum a_j \mu_j(y) + \sum a_\ell b_\ell \\ &= \omega(x,y) - \sum \lambda_k(x) \mu_k(y) = \tilde{\omega}(x,y).
\end{align*}

In the other direction, let $(X, Y, \omega, (\lambda_i), (\mu_i))$ be a nonmaximal complement datum compatible with a dual system for $\gl(V,V_*)$.  Fix $x \in X_0$ and $y \in Y_0$ such that $\tilde{\omega}(x,y)=1$.  Let $v_0 := x - \sum \lambda_j(x)v_j$ and $v^0 := y - \sum \mu_k(y)v^k$.  We calculate $\langle v_0, v^i \rangle = \lambda_i(x) - \lambda_i(x) = 0$ for $i>0$.  Similarly, $\langle v^i, v_0 \rangle = 0$ for $i>0$.  Finally, $\langle v_0, v^0 \rangle = \tilde{\omega}(x,y) = 1.$
\end{proof}

\begin{lemma} \label{lemma2}
Any nondegenerate complement datum is compatible with some (self-)dual system.
\end{lemma}
\begin{proof}
Let $(X, Y, \omega, (\lambda_i), (\mu_i))$ be a nondegenerate complement datum for $\gl_\infty$ or $\sl_\infty$.  Define $V := X \oplus \bigoplus_i \C v_i$ and $V_* := Y \oplus \bigoplus_i \C v^i$.  We define a bilinear map $\langle \cdot, \cdot \rangle$ on $V \times V_*$ by extending $\omega$ via $\langle x, v^i \rangle := \lambda_i(x)$ for $x \in X$, and $\langle v_i, y \rangle := \mu_i(y)$ for $y \in Y$, and $\langle v_i, v^j \rangle :=  \delta_{ij}$.  

Consider a nonzero element $v = x + \sum_i a_i v_i \in V$.  We see that $\langle v, v^j \rangle$ is nonzero for some $j$, unless $x \in X_0$ and a certain condition on $a_i$ is satisfied, namely that $a_i = -\lambda_i(x)$.  But if $v = x -\sum_i \lambda_i(x) v_i$ is nonzero, then $x$ must be nonzero, and there exists $y \in Y$ such that $0 \neq \tilde{\omega}(x,y) = \omega(x,y) -\sum_i \lambda{_i}(x) \mu_i(y) = \langle x , y \rangle - \sum_i \lambda_i(x) \langle v_i , y \rangle = \langle v, y \rangle$.  Similarly, any nonzero element in $V_*$ pairs nontrivially with an element of $V$.  Hence $\langle \cdot, \cdot \rangle$ is a nondegenerate pairing between $V$ and $V_*$.  So $\{\C v_i, \C v^i \}$ is a dual system for $\gl(V,V_*)$ or for $\sl(V,V_*)$ which is compatible with the given complement datum.

Now suppose we are given a nondegenerate complement datum $(X, \omega, (\lambda_i))$ for $\so_\infty$ or $\sp_\infty$.  Define $V := X \oplus \bigoplus_i \C v_i$.  We define a bilinear form $\langle \cdot, \cdot \rangle$ on $V$ by extending $\omega$ appropriately (i.e. symmetrically or antisymmetrically) via $\langle x, v_i \rangle := \lambda_i(x)$ for $x \in X$, and $\langle v_i, v_j \rangle := 0$ for $i \neq -j$, and $\langle v_i, v_{-i} \rangle := 1$ for $i > 0$.  

Consider a nonzero element $v = x + \sum_i a_i v_i \in V$.  We see that $\langle v, v_j \rangle$ is nonzero for some $j$, unless $x \in X_0$ and the coefficients $a_i$ satisfy some condition.  In the $\so_\infty$ case, the condition is $a_i = -\lambda_{-i}(x)$.  But if $v = x -\sum_i \lambda_{-i}(x) v_i$, then $x$ is nonzero, and there exists $y \in X$ such that $0 \neq \tilde{\omega}(x,y) = \omega(x,y) - \sum_i \lambda_{-i}(x) \lambda_i(y) = \langle x , y \rangle - \sum_i \lambda_{-i}(x) \langle v_i , y \rangle = \langle v, y \rangle$.  

In the $\sp_\infty$ case, the condition is $a_i = - \sgn(i) \lambda_{-i}(x)$.  But if $v = x -\sum_i \sgn(i) \lambda_{-i}(x) v_i$, then $x$ is nonzero, and there exists $y \in X$ such that $0 \neq \tilde{\omega}(x,y) = \omega(x,y) + \sum_i \sgn(i) \lambda_{-i}(x) \lambda_i(y) = \langle x , y \rangle - \sum_i \sgn(i) \lambda_{-i}(x) \langle v_i , y \rangle = \langle v, y \rangle$.  

Thus, in either case, $\langle \cdot, \cdot \rangle$ is a nondegenerate bilinear form.   So $\{\C v_i \}$ is a self-dual system for $\so(V)$ or $\sp(V)$ which is compatible with the given complement datum.
\end{proof}

\begin{defn}
\begin{enumerate}
\item \label{def1} Two complement data $(X, Y, \omega, (\lambda_i), (\mu_i))$ and $(X', Y', \omega', (\lambda'_i), (\mu'_i))$ for $\gl_\infty$ or $\sl_\infty$ are \emph{equivalent} if there exist a bijection $\sigma: \Z_{>0} \rightarrow \Z_{>0}$,  isomorphisms $\pi: X \rightarrow X'$ and $\pi: Y \rightarrow Y'$, and nonzero constants $\alpha_i$, such that $( \lambda_i (x) ) \approx ( \alpha_i \lambda'_{\sigma(i)} \circ \pi (x) )$ for all $x \in X$, $( \mu_i (y) ) \approx ( \alpha_i^{-1} \mu'_{\sigma(i)} \circ \pi (y) )$ for all $y \in Y$, and $\tilde{\omega} - \pi^*(\tilde{\omega}') = 0$ on $X \otimes Y$.  

\item \label{def2} Two complement data $(X, \omega, (\lambda_i))$ and $(X', \omega', (\lambda'_i))$ for $\so_\infty$ are \emph{equivalent} if there exist a bijection $\sigma : \Z_{\neq 0} \rightarrow \Z_{\neq 0}$ such that $\sigma (-i) = - \sigma (i)$, an isomorphism $\pi : X \rightarrow X'$, and nonzero constants $\alpha_i$ with $\alpha_{-i} = \alpha_i^{-1}$, such that $( \lambda_i (x) ) \approx ( \alpha_i \lambda'_{\sigma(i)} \circ \pi (x) )$ for all $x \in X$ and $\tilde{\omega} - \pi^*(\tilde{\omega}') = 0$ on $X \otimes X$.  

\item \label{def3} Two complement data $(X, \omega, (\lambda_i))$ and $(X', \omega', (\lambda'_i))$ for $\sp_\infty$ are \emph{equivalent} if there exist a bijection $\sigma : \Z_{\neq 0} \rightarrow \Z_{\neq 0}$ such that $\sigma (-i) = - \sigma (i)$, an isomorphism $\pi : X \rightarrow X'$, and nonzero constants $\alpha_i$ with $\alpha_{-i} = \alpha_i^{-1} \sgn(i) \sgn( \sigma(i))$, such that $( \lambda_i (x) ) \approx ( \alpha_i \lambda'_{\sigma(i)} \circ \pi (x) )$ for all $x \in X$ and $\tilde{\omega} - \pi^*(\tilde{\omega}') = 0$ on $X \otimes X$.  
\end{enumerate}
\end{defn}

Note that the left hand side of the equation $\tilde{\omega} - \pi^*(\tilde{\omega}') = 0$ in (\ref{def1}) above is in fact a finite sum when applied to any $x \otimes y \in X \otimes Y$.   Likewise, $\tilde{\omega} - \pi^*(\tilde{\omega}')$ in (\ref{def2}) and (\ref{def3}) is a finite sum when applied to any $x \otimes x' \in X \otimes X$.

The following theorem is our main result in Section \ref{sectionc}.

\begin{thm}\label{equivthm}
Let $\g$ be one of $\gl(V,V_*)$, $\sl(V,V_*)$, $\so(V)$, and $\sp(V)$.
The compatibility relation between (self-)dual systems and complement data induces a bijection between conjugacy classes of submaximal toral subalgebras of $\g$ and equivalence classes of nondegenerate complement data.   Under this bijection, conjugacy classes of maximal toral subalgebras of $\g$ correspond to equivalence classes of nondegenerate maximal complement data.
\end{thm}

\begin{proof}
We will prove the proposition in the case of $\g = \gl(V,V_*)$, and the other cases are similar.

There is a bijection between conjugacy classes of submaximal toral subalgebras of $\g$ and conjugacy classes of dual systems for $\g$, which comes from the bijection of submaximal toral subalgebras of $\g$ and dual systems for $\g$.  There is also a bijection between conjugacy classes of dual systems for $\g$ and equivalence classes of dual systems, by Lemma \ref{lemma3}.  We will show that compatibility induces a bijection between equivalence classes of dual systems and equivalence classes of nondegenerate complement data.  The second statement then follows immediately from Proposition \ref{compatibility}, which says that compatibility preserves maximality.

We first prove that compatibility induces a well-defined map from equivalence classes of dual systems to equivalence classes of complement data.  Let  $\{L_i, L^i \}$ be a dual system for $\gl(V,V_*)$, and $\{ M_i, M^i \}$ an equivalent dual system for $\gl(W, W_*)$.  Choose complement data $(X, Y, \omega, (\lambda_i), (\mu_i))$ and $(X', Y', \omega', (\lambda'_i), (\mu'_i))$ compatible with these dual systems.  We have implicitly chosen generators $v_i \in L_i$, $v^i \in L^i$, $w_i \in M_i$, and $w^i \in M^i$.  We have isomorphisms
 $\varphi: V \rightarrow W$ and $\varphi: V_* \rightarrow W_*$ and a bijection $\sigma: \Z_{>0} \rightarrow \Z_{>0}$ such that $\langle \cdot, \cdot \rangle = \varphi^*\langle \cdot, \cdot \rangle'$, $\varphi(L_i) =  M_{\sigma(i)}$, and $\varphi(L^i) = M^{\sigma(i)}$.  Let $\alpha_i \in \C$ be the nonzero constant such that $\varphi (v^i) = \alpha_i w^{\sigma(i)}$ (and hence $\varphi (v_i) = \alpha_i^{-1} w_{\sigma(i)}$).

The map
\begin{eqnarray*}
\varphi|_X : X & \rightarrow & X' \oplus \Span_i \{M_i \} \\
x & \mapsto & (\pi (x) , \sum \nu_i(x) w_i)
\end{eqnarray*}
defines a map $\pi: X \rightarrow X'$ and linear functionals $\nu_i:X \rightarrow \C$.  Notice that $( \nu_i (x) ) \approx 0$ for all $x \in X$.
Clearly $\pi$ is an isomorphism.
Similarly, 
\begin{eqnarray*}
\varphi|_Y : Y & \rightarrow & Y' \oplus \Span_i \{M^i\} \\
y & \mapsto & (\pi (y) , \sum \eta_i(y) w^i)
\end{eqnarray*}
defines $\pi : Y \rightarrow Y'$ and $\eta_i:Y \rightarrow \C$.  Again, $( \eta_i(y) ) \approx 0$ for all $y \in Y$.
Compute
\begin{eqnarray*}
\lambda_i (x) & = & \langle x , v^i \rangle  = \langle \varphi (x) , \varphi (v^i) \rangle' \\
& = & \langle \pi(x) + \sum \nu_j (x) w_j , \: \alpha_i w^{\sigma(i)} \rangle' 
 =  \alpha_i \lambda_{\sigma(i)}' (\pi (x)) + \alpha_i \nu_{\sigma(i)} (x).
\end{eqnarray*}
Hence for any $x \in X$ we have $( \lambda_i (x) ) \approx ( \alpha_i \lambda_{\sigma(i)}' \circ \pi (x) )$.  The analogous calculation shows $( \mu_i (y) ) \approx ( \alpha_i^{-1} \mu_{\sigma(i)}' \circ \pi (y) )$ for all $y \in Y$.
Then compute
\begin{eqnarray*}
\omega (x,y) & = & \langle x , y \rangle  =  \langle \varphi (x) , \varphi (y) \rangle' \\
& = & \langle \pi (x) + \sum \nu_j(x) w_j , \pi (y) + \sum \eta_k (y) w^k  \rangle'  \\
& = & \omega'( \pi (x) , \pi (y) ) - \sum \lambda_k'(\pi (x)) \mu_k'(\pi(y)) + \sum \lambda_l(x)\mu_l(y),
\end{eqnarray*}
i.e. $\tilde{\omega} - \pi^* \tilde{\omega}' = 0$ on $X \otimes Y$.  
This shows that the two complement data  $(X, Y, \omega, (\lambda_i), (\mu_i))$ and $(X', Y', \omega', (\lambda'_i), (\mu'_i))$ are equivalent.
Hence the compatibility relation induces a map from equivalence classes of dual systems to equivalence classes of nondegenerate complement data, where nondegeneracy follows from Proposition \ref{compatibility}.

The surjectivity of this map follows from Lemma \ref{lemma2}.  For injectivity we will show that if two dual systems are compatible with equivalent complement data, then the two dual systems are equivalent.

Suppose that the dual systems $\{ \C v_i , \C v^i \}$ for $\gl(V,V_*)$ and $\{ \C w_i , \C w^i \}$ for $\gl(W,W_*)$ are compatible with complement data $(X, Y, \omega, (\lambda_i), (\mu_i))$ and $(X', Y', \omega', (\lambda'_i), (\mu'_i))$, respectively, and that these complement data are equivalent.  There exist isomorphisms $\pi: X \rightarrow X'$ and $\pi: Y \rightarrow Y'$, nonzero constants $\alpha_i$ together with a bijection $\sigma: Z_{>0} \rightarrow Z_{>0}$ such that $( \lambda_i (x) ) \approx ( \alpha_i \lambda'_{\sigma(i)} \circ \pi (x) )$ for all $x \in X$ and $( \mu_i (y) ) \approx ( \alpha_i^{-1} \mu'_{\sigma(i)} \circ \pi (y) )$ for all $y \in Y$, and $\tilde{\omega} - \pi^*\tilde{\omega}' = 0$ on $X \otimes Y$.
Define linear functionals $\nu_i$ on $X$ and $\eta_i$ on $Y$ by
$$\nu_i: X \rightarrow \C , \: x \mapsto \lambda_i (x) - \alpha_i \lambda'_{\sigma(i)}(\pi(x))$$
and
$$\eta_i: Y \rightarrow \C , \: y \mapsto \mu_i (y) - \alpha_i^{-1} \mu'_{\sigma(i)}(\pi(y)).$$
Then $( \nu_i (x) ) \approx 0$ for all $x \in X$ and $( \eta_i (y) ) \approx 0$ for all $y \in Y$.
We have $V = X \oplus \bigoplus_i \C v_i$ and $V_* = Y \oplus \bigoplus_i \C v^i$.  Similarly, $W = X' \oplus \bigoplus_i \C w_i$ and $W_* = Y' \oplus \bigoplus_i \C w^i$.
We define
\begin{eqnarray*}
\varphi: V & \rightarrow & W \\
x + \sum a_i v_i & \mapsto & \pi(x) + \sum (\nu_i(x) + a_i ) \alpha_i^{-1} w_{\sigma(i)}
\end{eqnarray*}
and 
\begin{eqnarray*}
\varphi: V_* & \rightarrow & W_* \\
y + \sum a_i v^i & \mapsto & \pi(y) + \sum (\eta_i(y) + a_i ) \alpha_i w^{\sigma(i)}.
\end{eqnarray*}
These isomorphisms establish the equivalence of dual systems.
\end{proof}

We are ready now to draw some corollaries from Theorem \ref{equivthm}.  

In the spirit of \cite{N-P}, we introduce the five \emph{standard invariants} of an infinite dual system (and hence of an infinite-dimensional submaximal toral subalgebra or a Cartan subalgebra of $\gl_\infty$ or $\sl_\infty$): 
$$(\rank \langle \cdot , \cdot \rangle |_{T^\bot \times S^\bot},  \dim T^\bot, \dim S^\bot, \cod(T^\bot \oplus S), \cod(S^ \bot \oplus T) ),$$
where $S:= \bigoplus_i L_i$ and $T:=\bigoplus_i L^i$. 
The invariants take values in $\Z_{\geq 0} \cup \{ \aleph_0 \}$. If a dual system for $\gl(V,V_*)$ or $\sl(V,V_*)$ is compatible with the complement datum $(X, Y, \omega, (\lambda_i), (\mu_i))$, then one may check that the five standard invariants are
$$(\rank \tilde{\omega}|_{X_0 \times Y_0}, \dim X_0, \dim Y_0, \dim X/X_0, \dim Y/Y_0 ).$$  

Similarly, there are three  \emph{standard invariants} of an infinite self-dual system (and hence of an infinite-dimensional submaximal toral subalgebra or a Cartan subalgebra of $\so_\infty$ or $\sp_\infty$): 
$$(\rank \langle \cdot , \cdot \rangle |_{S^\bot \times S^\bot},  \dim S^\bot, \cod(S^ \bot \oplus S) ),$$
where $S:= \bigoplus_i L_i$.
If a self-dual system for $\so(V)$ or $\sp(V)$ is compatible with the complement datum $(X, \omega, (\lambda_i))$, then the three standard invariants are
$$(\rank \tilde{\omega}|_{X_0 \times X_0}, \dim X_0, \dim X/X_0 ).$$ 

\begin{cor} \label{sub}
\begin{enumerate}
\item There exists a submaximal toral subalgebra of $\gl_\infty$ or $\sl_\infty$ with standard invariants $(d,p,q,m,n)$ precisely when $0 \leq p-d \leq n$ and $0 \leq q-d \leq m$. 
\item There exists a submaximal toral subalgebra of $\so_\infty$ with standard invariants $(d,p,m)$ precisely when $0 \leq p-d \leq m$.
\item There exists a submaximal toral subalgebra of $\sp_\infty$ with standard invariants $(d,p,m)$ precisely when $0 \leq p-d \leq m$ and $d$ is even when it is finite.
\end{enumerate}
\end{cor}

\begin{proof}
We will give proofs which assume that the invariants are all finite, and we leave it to the reader to make the necessary modifications for the cases when the invariants are allowed to equal $\aleph_0$.
\begin{enumerate}
\item From the definition of a nondegenerate complement datum, $\tilde{\omega}$ yields surjections $X /X_0 \twoheadrightarrow (\mathrm{right \: ker \:} \tilde{\omega}|_{X_0 \times Y_0})^*$ and $Y/Y_0 \twoheadrightarrow ( \mathrm{left \: ker \:} \tilde{\omega}|_{X_0 \times Y_0})^*$.  Thus $\dim X / X_0 \geq \dim Y_0 - \rank \tilde{\omega}|_{X_0 \times Y_0} \geq 0$ and $\dim Y / Y_0 \geq \dim X_0 - \rank \tilde{\omega}|_{X_0 \times Y_0} \geq 0$.

Given $0 \leq p-d \leq n$ and $0 \leq q-d \leq m$, we construct a nondegenerate complement datum for $\gl_\infty$ or $\sl_\infty$ with standard invariants $(d,p,q,m,n)$ as follows.  Let $X$ be a vector space with basis $\{ x_{-p}, \ldots, x_{-1} , x_1 , \ldots, x_m \}$, and let $Y$ be a vector space with basis $\{ y_{-q}, \ldots, y_{-1} , y_1, \ldots, y_n \}$.  Define $\omega : X \times Y \rightarrow \C$ by setting $\omega(x_{-p+j},y_{-q+j}):=1$ for $j=0,\ldots,d-1$, and setting $\omega(x_{-j},y_j) := 1$ for $j=1,\ldots,p-d$, and setting $\omega(x_j,y_{-j}) :=1$ for $j=1,\ldots,q-d$ and letting all other pairings in this basis be trivial.  For $\iota \in \Z_{>0}$ define $\lambda_\iota := x_{\bar{\iota}}^*$ 
where $\bar{\iota} \in \{1, \ldots,m\}$ and $\bar{\iota} \equiv \iota (\mathrm{mod} \; m)$, and $\mu_\iota := y_{\bar{\iota}}^*$ 
where $\bar{\iota} \in \{1, \ldots,n\}$ and $\bar{\iota} \equiv \iota (\mathrm{mod} \; n)$.  
Then $X_0$ has basis $\{ x_{-p}, \ldots, x_{-1} \}$, and $Y_0$ has basis $\{ y_{-q}, \ldots, y_{-1} \}$, and $\tilde{\omega}|_{X_0 \times Y_0}$ has rank $d$.  This complement datum is nondegenerate, hence it gives rise to a submaximal toral subalgebra of $\gl_\infty$ or $\sl_\infty$ with standard invariants $(d,p,q,m,n)$.

\item From the definition of a nondegenerate complement datum, $\tilde{\omega}$ yields a surjection $X /X_0 \twoheadrightarrow (\rad \tilde{\omega}|_{X_0 \times X_0})^*$. Thus $\dim X / X_0 \geq \dim X_0 - \rank \tilde{\omega}|_{X_0 \times X_0} \geq 0$.

Given any $0 \leq p-d \leq m$, we construct a nondegenerate complement datum for $\so_\infty$ with standard invariants $(d,p,m)$ as follows.  Let $X$ be a vector space with basis $\{ x_{-p}, \ldots, x_{-1} , x_1 , \ldots, x_m \}$.  We define $\omega : X \times X \rightarrow \C$ by setting $\omega(x_{-p+j},x_{-p+j}):= 1$ for $j = 0, \ldots, d-1$ and setting $\omega(x_{-j},x_j) := 1$ for $j=1,\ldots,p-d$.  For $\iota \in \Z_{\neq 0}$ define $\lambda_\iota := x_{ \bar{\iota}}^*$ where $\bar{\iota} \in \{1, \ldots,m\}$ and $\bar{\iota} \equiv \iota (\mathrm{mod} \; m)$.  Then $X_0$ has basis $\{ x_{-p}, \ldots, x_{-1} \}$, and the rank of $\tilde{\omega}|_{X_0 \times X_0}$ is $d$.  This complement datum is nondegenerate, hence it gives rise to a submaximal toral subalgebra with standard invariants $(d,p,m)$.

\item The proof of the inequality is the same as for $\so_\infty$.  We see that $d$, being the rank of an antisymmetric bilinear form, is even if it is finite.

Given any $0 \leq p-d \leq m$ and $d$ even, we construct a nondegenerate complement datum for $\sp_\infty$ with standard invariants $(d,p,m)$ as follows.  Let $X$ be a vector space with basis $\{ x_{-p}, \ldots, x_{-1} , x_1 , \ldots, x_m \}$.  We define $\omega : X \times X \rightarrow \C$ by setting $\omega(x_{-p+2j},x_{-p+2j+1}):= 1$ for $j = 0, \ldots, d/2 -1$ and setting $\omega(x_{-j},x_j) := 1$ for $j=1,\ldots,p-d$.  For $\iota \in \Z_{\neq 0}$ define $\lambda_\iota := x_{\bar{\iota}}^*$ where $\bar{\iota} \in \{1, \ldots,m\}$ and $\bar{\iota} \equiv \iota (\mathrm{mod} \; m)$ .  Then $X_0$ has basis $\{ x_{-p}, \ldots, x_{-1} \}$, and the rank of $\tilde{\omega}|_{X_0 \times X_0}$ is $d$.  This complement datum is nondegenerate, hence it gives rise to a submaximal toral subalgebra with standard invariants $(d,p,m)$.

\end{enumerate}
\end{proof}

A submaximal toral subalgebra of $\gl_\infty$, $\sl_\infty$, or $\sp_\infty$ is maximal if and only if its first standard invariant is $0$.  A submaximal toral subalgebra of $\so_\infty$ is maximal if and only if its first standard invariant is $0$ or $1$.  Hence we immediately conclude the following corollary. 

\begin{cor}
\begin{enumerate}
\item \label{cora}There exists a Cartan subalgebra of $\gl_\infty$ or $\sl_\infty$ with standard invariants $(d,p,q,m,n)$ precisely when $d=0$ and $p \leq n$ and $q \leq m$. 

\item There exists a Cartan subalgebra of $\so_\infty$ with standard invariants $(d,p,m)$ precisely when $d=0,1$ and $d \leq p \leq m+d$.

\item There exists a Cartan subalgebra of $\sp_\infty$ with standard invariants $(d,p,m)$ precisely when $d=0$ and $p \leq m$.
\end{enumerate}
\end{cor}

\begin{cor}
\begin{enumerate}
\item \label{corb} The cardinality of the set of conjugacy classes of Cartan subalgebras of $\gl_\infty$ or $\sl_\infty$ with given standard invariants is $0$, $1$, $2$, or continuum.  One conjugacy class occurs precisely for standard invariants $(0,0,0,0,0)$.  Two conjugacy classes occur precisely for standard invariants $(0,0,0,1,0)$, $(0,0,0,0,1)$, $(0,1,0,0,1)$, and $(0,0,1,1,0)$.

\item The cardinality of the set of conjugacy classes of Cartan subalgebras of $\so_\infty$ with given standard invariants is $0$, $1$, or continuum.  One conjugacy class occurs precisely for standard invariants $(0,0,0)$ and $(1,1,0)$.
 
\item The cardinality of the set of conjugacy classes of Cartan subalgebras of $\sp_\infty$ with given standard invariants is $0$, $1$, or continuum.   One conjugacy class occurs precisely for standard invariants $(0,0,0)$.
\end{enumerate}
\end{cor}

\begin{proof}
\begin{enumerate}
\item  Let $(X,Y,\omega, (\lambda_i),(\mu_i))$ be a complement datum for  $\gl_\infty$ or $\sl_\infty$  with standard invariants $(0,p,q,m,n)$. After possibly interchanging the role of $X$ and $Y$, we may assume that $m \geq n$.  First we consider each of the special cases for $\gl_\infty$, and then prove that in all other nontrivial cases the set of conjugacy classes has cardinality continuum.  Let  $\{e_i \}$ and $\{e^i \}$ be dual bases of $V$ and $V_*$, indexed by $\Z_{>0}$.

\begin{enumerate} 
\item[$(0,0,0,0,0)$] 

Here both $X$ and $Y$ have dimension $0$, so there is exactly one linear functional on each space and exactly one pairing between them.  Thus there is exactly one complement datum, and it is nondegenerate and maximal.  Up to equivalence there is still exactly one complement datum.  A dual system in this class is $\{ L_i = \C e_i , L^i = \C e^i : i \in \Z_{>0} \}$.

\item[$(0,0,0,1,0)$]

We see that $X$ is $1$-dimensional and $Y$ is $0$-dimensional.  By setting
$$(\lambda_i , \mu_i) := 
 \begin{cases}
(1,0) & i \in \Z_{>0} \\
(\bar{i},0) & i \in \Z_{>0},
\end{cases}
$$
where $\bar{i} \equiv i(\mathrm{mod} \; 2)$ and $\omega = 0$, we define two complement data.  These are inequivalent because whether $\lambda_i$ is zero infinitely often or finitely often is an invariant.

Notice that in the case under consideration $\mu_i=0$ and $\omega =0$ for an arbitrary nondegenerate maximal complement datum $(X,Y, \omega, (\lambda_i) , (\mu_i))$.  Up to equivalence, we may assume that $\lambda_i \in \{0,1\}$, where $\lambda_i =1$ must occur infinitely often since $X_0 = 0$.  Therefore, up to equivalence, there are just two possibilities: $\lambda_i=0$ occurs infinitely often or finitely often.  These equivalence clases are represented by the above complement data.  A dual system in the first class is $\{ L_i = \C (e_i - e_{i+1}) , L^i = \C\left( \sum_{k=1}^{i} e^k \right) : i \in \Z_{>0} \}$, and a dual system in the second class is $\{ L_i = \C(e_i - \bar{i} e_{i+2}), L^i = \C \left(e^{i} + \bar{i} \sum_{k=1}^{(i-1)/2} e^{2k-1}\right) : i \in Z_{>0} \}$, where $\bar{i} \in \{0,1\}$ and $\bar{i} \equiv i (\mathrm{mod} \; 2)$.

\item[$(0,0,1,1,0)$]

We see that $X$ is $1$-dimensional and $Y = Y_0$ is $1$-dimensional.  We define two complement data by setting
$$(\lambda_i , \mu_i) := 
 \begin{cases}
(1,0) & i \in \Z_{>0} \\
(\bar{i},0) & i \in \Z_{>0},
\end{cases}
$$
where $\bar{i} \equiv i(\mathrm{mod} \; 2)$ and $\omega = 1$.  These are inequivalent because whether or not $\lambda_i = 0$ occurs infinitely often is an invariant.

Since $Y = Y_0$, $( \mu_i ) \approx 0$, so every complement datum is equivalent to one with $\mu_i = 0$ for all $i$.  Fix $x \in X$, and we may assume $\lambda_i(x) \in \{0,1\}$ for all $i$.  Since $X_0 = 0$, we must have $\lambda_i(x) = 1$ for infinitely many $i$.  Now $\tilde{\omega} = \omega$ is nonzero by nondegeneracy, so we may fix $y \in Y$ such that $\tilde{\omega}(x,y) =1$.  Up to equivalence, there are just two possibilities:  $\lambda_i (x)= 0$ occurs infinitely often or finitely often.  A dual system in the first class is $\{ L_i = \C (e_{i+1}), L^i = \C ( e^1 + e^{i+1}) : i \in \Z_{>0} \}$, and a representative of the second class is $\{ L_i = \C(e_{i+1}) , L^i = \C (\bar{i}e^1 + e^{i+1}) : i \in \Z_{>0} \}$, where $\bar{i} \in \{0,1\}$ and $\bar{i} \equiv i (\mathrm{mod} \; 2)$
\end{enumerate}

Now suppose that $p \leq n$ and $q \leq m$ and that we are not in one of the above cases.  Thus either  $m=n=1$ or $m \geq 2$.  Since the cardinality of the set of all complement data for $X$ and $Y$ is at most continuum, it is enough to show that the cardinality of equivalence classes in each of these cases is at least continuum.

\begin{enumerate}

\item[$m=n=1$]

We will construct a family of complement data which represents continuum many equivalence classes.  The space $(X/X_0)^* \otimes (Y/Y_0)^*$ is a vector space of dimension $1$ which we identify with $\C$.  For any $z \in \C \setminus \{0, 1\}$ we define a complement datum $D_z$ as follows.  Choose $\lambda_i$ and $\mu_i$ to vanish on $X_0$ and $Y_0$ respectively, such that $\lambda_i \otimes \mu_i$ is $1$ or $z$ and each of these two values occurs infinitely often.  Choose $\omega$ to give surjections $X/X_0 \twoheadrightarrow Y_0^*$ and $Y/Y_0 \twoheadrightarrow X_0^*$.  Since $\tilde{\omega}$ coincides with $\omega$ when restricted to $X \otimes Y_0$ or $X_0 \otimes Y$,  such a complement datum $D_z$ is nondegenerate and maximal for any $z$.

For $D_z$ and $D_{z'}$ to be equivalent, there must be an element $\pi \in \GL(\C)$ and a permutation $\sigma$ such that $(\pi(\lambda_i\otimes \mu_i)) \approx (\lambda'_{\sigma(i)} \otimes \mu'_{\sigma(i)})$.  Thus there must be an element of $\C^\times$ sending $\{1, z\}$ to $\{1, z'\}$.  Hence the only other element from this family in the equivalence class of $D_z$ is $D_{1/z}$. Therefore there are continuum many equivalence classes.  

\item[$m \geq 2$] 

Consider the projective space $P = \Proj((X/X_0)^*)$, and let $S \subset P$ be an at most countable subset. We define a complement datum $D_S$ as follows.  Choose $0 \neq \lambda_i \in X^*$ such that for each $s \in S$, $\lambda_i$ lies on the line $s$ for infinitely many $i$.  Choose $\mu_i \in Y^*$ such that that $(\mu_i (y)) \approx 0$ implies $y \in Y_0$.   Finally, choose $\omega$ to give surjections $X/X_0 \twoheadrightarrow Y_0^*$ and $Y/Y_0 \twoheadrightarrow X_0^*$.  Since $\tilde{\omega}$ coincides with $\omega$ when restricted to $X \otimes Y_0$ or $X_0 \otimes Y$,  such a complement datum $D_S$ is nondegenerate and maximal for any $S$.  Clearly if $D_S$ and $D_{S'}$ are equivalent, then there must be an element $\pi \in \GL(X)$, a bijection $\sigma$, and nonzero constants $\alpha_i$ such that $(\lambda_i) \approx (\alpha_i \lambda'_{\sigma(i)} \circ \pi)$.  Thus there must be an element of $\Aut(P)$ mapping $S$ to $S'$.  Hence it is enough to show that there are uncountably many equivalence classes under $\Aut(P)$ of at most countable subsets of $P$.

Suppose that $m$ is finite.  Fix a set of $m+1$ points in general position in $P$.  For each $z \in P$ in general position, let $S_z$ be the union of the fixed set and $\{z\}$.  Since $\Aut(\Proj^{m-1})$ acts with finite stabilizers on sets of $m+1$ points in $\Proj^{m-1}$, each $S_z$ is equivalent to at most a finite number of other such $S_{z'}$.  Thus there are at least continuum many equivalence classes.

Suppose that $m > 1$, where we allow $m = \aleph_0$.  The following argument was suggested by Scott Carnahan and Anton Geraschenko.  Choose countably many distinct lines $\ell_1, \ell_2, \ldots$ in $P$, which is possible since $m > 1$.  Fix a subset $T$ of $\Z_{>2}$.  Let $t_i$ be the $i$th smallest element of $T$.  On $\ell_1$ choose $t_1$ distinct points.  On $\ell_2$ choose $t_2$ distinct points.  On $\ell_3$ choose $t_3$ distinct points each of which is not colinear with any pair of the chosen points on $\ell_1$ and $\ell_2$.  This is possible because there are only finitely many lines going through two of the previously chosen points.  Continue choosing points to get a set $S_T$ in $P$ satisfying the property that exactly one subset of $t_i$ points lie on the same line in $P$, and such that $t_i$ points lie on $\ell_i$.  Since linear transformations preserve colinearity, we see that $S_T$ is equivalent to $S_{T'}$ if and only if $T = T'$.  Since there are continuum many subsets of $\Z_{>2}$, we have constructed representatives of continuum many equivalence classes.

\end{enumerate}

\item Let $(X,\omega,(\lambda_i))$ be a complement datum for $\so_\infty$ with standard invariants $(d,p,m)$.  There are two special cases for $\so_\infty$, and then we prove that in all other nontrivial cases, there are continuum many conjugacy classes.

\begin{enumerate}

\item[$(0,0,0)$] 
In this case $X_0 = X = 0$, hence there is exactly one equivalence class of complement data.  A Cartan subalgebra with these invariants arises naturally from the realization of $\so_\infty$ as $\bigcup \so_{2n}$.

\item[$(1,1,0)$] 
In this case $X_0 = X$ is $1$-dimensional and $\tilde{\omega}$ is nontrivial.  Since $(\lambda_i) \approx 0$, there is exactly one equivalence class of complement data.  A Cartan subalgebra with these invariants arises naturally from the realization of $\so_\infty$ as $\bigcup \so_{2n+1}$.

\item[$m \geq 1$] 
Consider the vector space $U = (X/X_0)^* \otimes (X/X_0)^*$.  Consider an at most countable set $S \subset U$. We define a complement datum $D_S$ as follows.  Choose $\lambda_i$ to vanish on $X_0$ such that $\lambda_i \otimes \lambda_{-i} \in S$ and each element of $S$ occurs for infinitely many $i$.  Choose $\omega$ to have rank $d$ on $X_0$ and to give a surjection $X /X_0 \twoheadrightarrow (\rad \omega |_{X_0 \times X_0})^*$.  Since $\tilde{\omega}$ coincides with $\omega$ when restricted to $X_0 \otimes X$,  such a complement datum $D_S$ is nondegenerate and maximal for any $S$. Clearly if $D_S$ and $D_{S'}$ are equivalent, then there must be an element $\pi \in \GL(X)$ and a bijection $\sigma$ such that $(\lambda_i \otimes \lambda_{-i}) \approx (\lambda'_{\sigma(i)} \circ \pi \otimes \lambda'_{-\sigma(i)} \circ \pi)$.  Thus there must be an element of $\Aut(U)$ mapping $S$ to $S'$.  Hence it is enough to show that there are uncountably many equivalence classes under $\Aut(U)$ of at most countable subsets of $U$.

Now we can use the same two arguments from the analogous case for $\gl_\infty$.  If $m$ is finite, we can choose $m^2+1$ points in general position and proceed as before.  If $m>1$, including $m = \aleph_0$, then we can choose points on lines in $U$ satisfying certain collinearity properties as before.
\end{enumerate}

\item  Let $(X,\omega,(\lambda_i))$ be a complement datum for $\sp_\infty$ with standard invariants $(0,p,m)$. There is one special case for $\sp_\infty$.

\begin{enumerate}

\item[$(0,0,0)$] In this case $X_0 = X = 0$, thus there is exactly one conjugacy class of Cartan subalgebras, and they arise naturally from the realization of $\sp_\infty$ as $\bigcup \sp_{2n}$.

\item[$m \geq 1$] Repeat the proof of the $m \geq 1$ case for $\so_\infty$, choosing $\lambda_i$ such that $\sgn(i) \lambda_i \otimes \lambda_{-i} \in S$. 

\end{enumerate}
\end{enumerate}
\end{proof}

In particular, for $\so_\infty$ and $\sp_\infty$, if there are finitely many conjugacy classes of Cartan subalgebras with given standard invariants, then the Cartan subalgebras with these standard invariants are splitting.

We conclude the paper with a description of a case with uncountably many equivalence classes.  Let  $\{e_i \}$ and $\{e^i \}$ be dual bases of $V$ and $V_*$, indexed by $\Z_{>0}$.  Consider for any pair of sequences $(a_i)$ and $(b_i)$ the dual system for $\gl(V,V_*)$ given by $\{ L_i = \C( b_i e_1+e_{i+2} ), L^i = \C( a_i e^2 +e^{i+2} ) : i \in \Z_{>0} \}$.  This dual system has invariants $(0,1,1,1,1)$ if and only if the sequences $(a_i)$ and $(b_i)$ are not almost equal to $0$.   The following are three binary invariants of these dual systems: whether for infinitely many $i$, $a_i$ is zero and $b_i$ is nonzero; whether for infinitely many $i$, $a_i$ is nonzero and $b_i$ is zero; and whether for infinitely many $i$, both $a_i$ and  $b_i$ are zero.  

We define a \emph{multiset} on a set $S$ as a map $m : S \rightarrow \Z_{\geq 0}$ such that $m$ vanishes outside a countable subset of $S$.  We say $m(x)$ is the \emph{multiplicity} of $x$.  Suppose the group $\C^\times$ acts on $S$.  Two multisets $m$, $m'$ on $S$ are \emph{almost proportional} if there exists $c \in C^\times$ such that for all $x \in S$, $| m(x) - m'(cx) | < \infty$, and $m(x) = m'(cx)$ for all but finitely many $x \in S$.  A sequence $(\alpha_i \in S)$ gives rise to a multiset $m : S \rightarrow \Z_{\geq 0}$ by setting $m(x) := |\{ i : \alpha_i = x \}|$.  Two sequences are  \emph{almost proportional} if the multisets arising from them are almost proportional.  Consequently, two sequences $(\alpha_i)$ and $(\beta_i)$ almost proportional if and only if there exists $c \in \C^\times$ and a permutation $\sigma$ such that $( \alpha_i ) \approx ( c \beta_{\sigma(i)} )$.

\begin{prop}
Every dual system with standard invariants $(0,1,1,1,1)$ is equivalent to one of the form $\{ L_i = \C( b_i e_1+e_{i+2} ), L^i = \C( a_i e^2 +e^{i+2} ) : i \in \Z_{>0} \}$.  Two such dual systems are conjugate if and only if their three binary invariants agree and the sequence $(a_i  b_i)$ is almost proportional to the sequence $(a_i'  b_i' )$.
\end{prop}

\begin{proof}
Fix a dual system with standard invariants $(0,1,1,1,1)$, and choose a compatible complement datum $(X', Y', \omega', (\lambda'_i), (\mu'_i))$ which satisfies the further condition that the $\lambda'_i$ vanish on $X'_0$ and the $\mu'_i$ vanish on $Y'_0$.  Choose nonzero vectors $x_0 \in X'_0$ and $y_0 \in Y'_0$.  Since the complement datum is nondegenerate, we can choose vectors $x \in X'$ and $y \in Y'$ such that $\omega'(x, y_0) = \omega'(x_0, y) = 1$.  As the complement datum is maximal, we have $\omega'(x_0, y_0) = 0$.  After possibly replacing $x$ with $x + \beta x_0$, we may assume without loss of generality that $\omega'(x,y) = 0$.  Let $a_i := \lambda'_i(x)$ and $b_i := \mu'_i(y)$.  

Consider the dual system $\{ L_i = \C v_i, L^i = \C v^i: i \in \Z_{>0} \}$, where $v_i := b_i e_1+e_{i+2}$ and $v^i := a_i e^2 +e^{i+2}$.  Let $X:= \Span \{e_1,e_2\}$ and $Y:= \Span \{e^1,e^2\}$.  Define a map $\pi : X' \rightarrow X$ by $x_0 \mapsto e_1$, $x \mapsto e_2$, and define $\pi : Y' \rightarrow Y$ by $y_0 \mapsto e^2$, $y \mapsto e^1$.  It is clear that $\pi$ is an equivalence of dual systems.  Thus every equivalence class of dual systems has a representative of the required form.

It remains to show that two such dual systems are equivalent precisely when they have the same binary invariants and the sequences $(a_i  b_i)$ and $(a_i'  b_i' )$ are almost proportional.  Suppose two dual systems, given by sequences $(a_i)$, $(b_i)$, $(a'_i)$, and $(b'_i)$, are equivalent.  
There exist a bijection $\sigma: \Z_{>0} \rightarrow \Z_{>0}$, isomorphisms $\pi: X \rightarrow X$ and $\pi: Y \rightarrow Y$, and nonzero constants $\alpha_i$, such that $( \lambda_i ) \approx ( \alpha_i \lambda'_{\sigma(i)} \circ \pi )$ and $( \mu_i ) \approx ( \alpha_i^{-1} \mu'_{\sigma(i)} \circ \pi )$.  We adopt the notation $x_0 := e_1$, $x := e_2$, $y_0 := e^2$, and $y := e^1$.  Suppose $\pi : X \rightarrow X$ is given by $x_0 \mapsto a x_0 + b x$, $x \mapsto c x_0 + d x$.  Then $(\lambda_i (x_0) ) \approx ( \alpha_i \lambda'_{\sigma(i)} (a x_0 + b x))$ implies $0 \approx ( b \alpha_i a'_i )$, hence $b=0$.  Since $\pi$ is an isomorphism, it follows that $a$ and $d$ are nonzero.  We also have $( \lambda_i(x) ) \approx ( \alpha_i \lambda_{\sigma(i)}' (c x_0 + d x) )$, thus $( a_i ) \approx ( \alpha_i d a_{\sigma(i)}' )$.  Similarly, there is a nonzero $\delta$ such that $( b_i ) \approx ( \alpha_i^{-1} \delta b_{\sigma(i)}' )$.  Therefore the binary 
 invariants agree.  Moreover $(a_i b_i ) \approx ( d \delta a_{\sigma(i)}' b_{\sigma(i)}' )$, so the sequences $(a_i  b_i)$ and $(a_i'  b_i' )$ are almost proportional.

Conversely, suppose that $(a_i b_i ) \approx ( c a_{\sigma(i)}' b_{\sigma(i)}' )$ with $c$ nonzero and that the binary invariants agree.  Choose nonzero constants $\alpha_i$ such that $( a_i ) \approx (  \alpha_i c a_{\sigma(i)}' )$ and $( b_i ) \approx ( \alpha_i^{-1} b_{\sigma(i)}'  )$.  This will be possible precisely because the sequences' binary invariants agree.  Let $\gamma$ be the constant $ \gamma := c^{-1} \sum ( c  a_{\sigma(i)}' b_{\sigma(i)}' - a_i b_i) $.  Define $\pi : X \rightarrow X$ by $x_0 \mapsto x_0$ and $x \mapsto cx$, and define $\pi : Y \rightarrow Y$ by $y_0 \mapsto c^{-1} y_0$ and $y \mapsto \gamma y_0 + y$.  One calculates that $( \lambda_i ) \approx ( \alpha_i \lambda'_{\sigma(i)} \circ \pi )$, $( \mu_i ) \approx ( \alpha_i^{-1} \mu'_{\sigma(i)} \circ \pi )$, and $\tilde{\omega} - \pi^*(\tilde{\omega}') = 0$.  
\end{proof}

\begin{figure}
\begin{center}
\caption{Parameter space for $(a_i,b_i)$}
\includegraphics{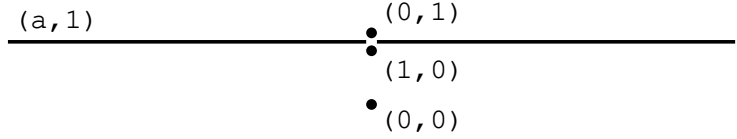} 
\label{fig2}
\end{center}
\end{figure}

We can visualize the preceding proposition as follows.  Cartan subalgebras of $\gl_\infty$ with standard invariants $(0,1,1,1,1)$  are in bijection with certain "admissible" multisets on the parameter space depicted in Figure \ref{fig2}, modulo almost proportionality.  More precisely, the parameter space is $(\C \times \C) / \C^\times$, where the action of $\C^\times$ on $\C \times \C$ is $c \cdot (x,y) = (cx,c^{-1}y)$.  A multiset is admissible if infinitely many points lie on each of the two lines.  The action of $\C^\times$ on the parameter space is $c \cdot [(x,y)] = [(cx,cy)] =[(c^2x,y)]$, and hence elements of $\C^\times$ simply rescale the line.

A similar analysis holds more generally.  Let $S$ be the quotient $((X/X_0)^* \times (Y/Y_0)^*) / \C^\times$, where $\C^\times$ acts on $(X/X_0)^* \times (Y/Y_0)^*$ by $c \cdot (x,y) = (cx,c^{-1}y)$.  There is a residual diagonal action of $\C^\times$ on $S$, which enables us to define almost proportional multisets on $S$.  If two complement data are equivalent, then the images of $(\lambda_i, \mu_i)$ and $(\lambda'_i, \mu'_i)$ in $S$ must be almost proportional.  However, the converse does not hold, as different choices of $\omega$ often yield inequivalent dual systems.  It is difficult to clarify this dependence in general;  nonetheless, in practice one should be able to distinguish nonconjugate dual systems.  For example, consider complement data with standard invariants $(0,0,0,1,1)$ and a fixed sequence $(\lambda_i, \mu_i)$ of $S$ with no unusual symmetries.  It is easy to see that any two choices of $\omega$ yield inequivalent dual systems.


\begin{thebibliography}{}
\bibitem[\textbf{BB}]{BahturinBenkart}Y.  Bahturin, G. Benkart,  {\it Highest weight modules for locally finite Lie algebras}, AMS/IP Stud. in Adv. Math. 4 (1997), 1--31.
\bibitem[\textbf{BS}]{BahturinStrade} Y. Bahturin, H. Strade,  { \it Some examples of locally finite simple Lie algebras}, Arch. Math. (Basel) 65 (1995), 23--26.
\bibitem[\textbf{B1}]{B1} A. Baranov, {\it Diagonal locally finite Lie algebras and a version of Ado's theorem}, J. Algebra 199 (1998), 1--39. 
\bibitem[\textbf{B2}]{B2} A. Baranov, {\it Simple diagonal locally finite Lie algebras}, Proc. London Math. Soc. 77 (1998), 362--386.
\bibitem[\textbf{BZ}]{BaranovZhilinski} A. Baranov, A. Zhilinski, {\it Diagonal direct limits of simple Lie algebras}, Comm. Algebra 27 (1998), 2249--2766. 
\bibitem[\textbf{B}]{Bourbaki} N. Bourbaki, {\it Groupes et Alg\`ebres de Lie}, Hermann, Paris, 1975.
\bibitem[\textbf{DP1}]{D-P1} I. Dimitrov, I. Penkov, {\it Weight modules of direct limit Lie algebras}, Intern. Math. Res. Notices 1999, No. 5, 223-249, math.RT/9808045.
\bibitem[\textbf{DP2}]{D-P2} I. Dimitrov, I. Penkov, {\it Borel subalgebras of $\gl(\infty)$}, Resenhas IME-USP 6 (2004), No. 2/3, 153-163.
\bibitem[\textbf{M}]{Mackey} G. Mackey, {\it On infinite dimensional linear spaces}, Trans. AMS 1945, Vol. 57, No. 2, 155-207.
\bibitem[\textbf{NP}]{N-P} K.-H. Neeb, I. Penkov, {\it Cartan subalgebras of $\gl_\infty$}, Canad. Math. Bull. 46(2003), 597-616.
\end{thebibliography}
\end{document}